\documentclass[11pt]{amsart}
\usepackage{amsmath,amssymb,amsfonts,url,mathptmx}
\usepackage[mathscr]{euscript}
 \usepackage{hyperref}
 \usepackage{mathrsfs}
\usepackage{graphicx}
\usepackage{color}
  \usepackage{epsfig}
\numberwithin{equation}{section}

\newtheorem{thm}{Theorem}[section]
\newtheorem{lem}{Lemma}[section]
\newtheorem{rem}{Remark}[section]
\newtheorem{prop}{Proposition}[section]

\newcommand{\hdot}{^\text{\r{}}\hspace{-.33cm}H}

\begin{document}
\title[Singular Liouville system]{Structure of bubbling solutions of Liouville systems with negative singular sources } 

\author{Yi Gu}
\address{Research Institute of Nuclear Power Operation (Shanghai)\\
        5th Floor, Building A, 1177 Shibo Avenue\\
        Pudong, Shanghai, China}
\email{guyi168@126.com}

\author{Lei Zhang}\footnote{Lei Zhang is partially supported by a Simons Foundation Collaboration Grant}
\address{Department of Mathematics\\
        University of Florida\\
        358 Little Hall P.O.Box 118105\\
        Gainesville FL 32611-8105}
\email{leizhang@ufl.edu}

\date{\today}

\begin{abstract}
 Liouville systems on Riemann surfaces are instrumental in modeling species growth and particle dynamics in biology and physics. Previously, we established a priori estimates for parameters across regions defined by critical hyper-surfaces. Here, we extend this by giving a priori estimates when parameters are critically positioned. This involves thoroughly characterizing bubble interaction, a key challenge in Liouville systems. During blowup events, we ascertain the exact heights of bubbling solutions about each blowup point, the integrals of each component, and the blowup points' positions. Moreover, as the parameter $\rho$ approaches a critical hyper-surface, we identify a pivotal leading term vital for numerous applications.\end{abstract}

\maketitle

\numberwithin{equation}{section}
\allowdisplaybreaks

\section{Introduction}

In this article we consider the following Liouville system defined on Riemann surface $(M,g)$:
\begin{equation}\label{1.1}
\Delta_g v_i+\sum_{j=1}^n a_{ij}\rho_j(\frac{H_je^{v_j}}{\int_M H_je^{v_j}dV_g}-1)=\sum_{l=1}^{N}4\pi\gamma_l(\delta_{p_l}-1),\quad i=1,..,n.
\end{equation}
where $H_1,\cdots,H_n$ are positive smooth functions on $M$, $p_1,...,p_N$ are distinct points on $M$, $4\pi \gamma_i \delta_{p_i}$ ($i=1,...,N$) are Dirac masses placed at $p_i$ with each $-1<\gamma_i<0$, $\rho_1,\cdots,\rho_n$ are nonnegative constants and without loss of generality, we assume $Vol(M)=1$. $\Delta_g$ is the Laplace-Beltrami operator ($-\Delta_g\geq 0$). Equation (1.1) is called Liouville system if all the entries in the coefficient matrix $A=(a_{ij})_{n\times n}$ are nonnegative. Here we point out that we assume the singular source on the right hand side is the same for all $i$ for simplicity.

 Liouville systems have significant applications across various fields. In geometry, when the system reduces to a single equation ($n=1$), it generalizes the renowned Nirenberg problem, which has been extensively researched over the past few decades (see \cite{barto2,barto3,bart-taran-jde,bart-taran-jde-2,ccl,chen-li-duke,kuo-lin-jdg,li-cmp,li-shafrir,wei-zhang-adv,wei-zhang-plms,zhangcmp,zhangccm}). In physics, Liouville systems emerge from the mean field limit of point vortices in the Euler flow (see \cite{biler,wolansky1,wolansky2,yang}) and are intricately linked to self-dual condensate solutions of the Abelian Chern-Simons model with $N$ Higgs particles \cite{kimleelee,Wil}. In biology, they appear in the stationary solutions of the multi-species Patlak-Keller-Segel system \cite{wolansky3} and are important for studying chemotaxis \cite{childress}. Understanding bubbling solutions, a significant challenge within Liouville systems, is crucial for advancing related fields.

The equation (\ref{1.1}) is usually written in an equivalent form by removing the singular sources on the right hand side. Let $G(x,y)$ be the Green's function of $-\Delta_g$ on $M$:
$$-\Delta_g G(x,p)=\delta_p-1, \quad \int_M G(x,p)dV_g(x)=0, $$
Note that in a neighborhood of $p$, 
$$G(x,p)=-\frac 1{2\pi}\log |x-p|+\gamma(x,p).$$
Using
$$u_i=v_i+4\pi\sum_{l=1}^{L_2}\gamma_lG(x,p_l),\quad i\in I:=\{1,...,n\}$$
we can write (\ref{1.1}) as
\begin{equation}\label{main-sys-2}
\Delta_gu_i+\sum_{j=1}^na_{ij}\rho_j(\frac{\mathfrak{h}_je^{u_j}}{\int_M\mathfrak{h}_je^{u_j}}-1)=0,\quad i\in I:=\{1,...,n\},
\end{equation}
where
$$\mathfrak{h}_i(x)=H_i(x)exp\{-\sum_{l=1}^{L_2}4\pi \gamma_lG(x,p_l)\}. $$
In a neighborhood around each singular source, 
say, $p_l$, in local coordinates, $\mathfrak{h}_j$ can be written as
$$\mathfrak{h}_j(x)=|x|^{2\gamma_l}h_j(x)$$
for some positive, smooth function $h_j(x)$.

Let $u=(u_1,...u_n)$ be a solution of (\ref{main-sys-2}), then it is standard to say $u$ belongs to  $\hdot^{1,n}(M)=\hdot^1(M)\times....\times \hdot^1(M)$ where 
$$ \hdot^{1}(M):=\{v\in L^2(M);\quad \nabla v\in L^2(M), \mbox{and }\,\, \int_M v dV_g=0\}. $$
Corresponding to \, $\hdot^{1,n}$ there is a variational form whose Euler-Lagrange equation is (\ref{main-sys-2}).

In \cite{linzhang1,linzhang2}, Lin and the second author completed a degree counting program to regular Liouville systems ( no singular source) under the following two assumptions on the
matrix $A$: (recall that $I=\{1,...,n\}$)
\begin{align*}
&(H1): \quad A\mbox{ is symmetric, nonnegative, irreducible and invertible.}\\
&(H2): \quad a^{ii}\leq 0,\,\, \forall i\in I,\quad  a^{ij}\geq 0 \,\, \forall i\neq j\in I,  \quad  \sum_{j=1}^na^{ij}\geq 0 \,\,\forall i\in I,
\end{align*}
where $(a^{ij})_{n\times n}$ is $A^{-1}$.  $(H1)$ is a rather standard assumption for Liouville systems, $(H2)$ says the interaction between equations is strong, for example when $n=2$, the non-negative matrix
$\left(\begin{array}{cc}
a_{11} & a_{12} \\
a_{12} & a_{22}
\end{array}
\right)$ satisfies $(H1)$ and $(H2)$ if  $a_{12}>\max\{a_{11},a_{22}\}$.
Later Lin-Zhang's work has been extended by the authors \cite{gu-zhang} for singular Liouville systems. Among other things we prove the following: Let $\Sigma$ be a set of critical values:
$$\Sigma:=\{8m\pi+\sum_{p_l\in \Lambda}8\pi \mu_l;  \Lambda\subset\{p_1,\cdots,p_N\}, m\in \mathbb{N}^+\cup\{0\}\}\setminus \{0\},  $$
where $\mu_l=1+\gamma_l$ and $\mathbb N^+$ is the set of natural numbers. If $\Sigma$ is written as
$$\Sigma=\{8\pi n_1<8\pi n_2<\cdots <8\pi n_L<...\},$$
then for $\rho=(\rho_1,...,\rho_n)$ satisfying
\begin{equation}\label{rho}
8\pi n_L \sum_{i\in I} \rho_i<\sum_{i,j\in I}a_{ij}\rho_i\rho_j<8\pi n_{L+1} \sum_{i\in I}\rho_i,
\end{equation}
there is a priori estimate for all solution $u$ to (\ref{1.1}), and the Leray-Schaudar degree $d_\rho$ for equation (\ref{1.1}) is computed. The main result in \cite{gu-zhang} states that if the parameter $\rho=(\rho_1,...,\rho_n)$ is not on any of the critical hyper-surface $\Gamma_L$ below and if the manifold has a non-positive Euler characteristic, then there exists a solution. 
Let
\begin{equation*}
\Gamma_{L}=
\left\{
\rho\,\big|\, \rho_i>0, i\in I;\,\,  \Lambda_L(\rho)=0,
\right\},
\end{equation*}
be the L-th hypersurface, where 
$$\Lambda_L(\rho)=4\sum_{i=1}^n\frac{\rho_i}{2\pi n_L}-\sum_{i=1}^n\sum_{j=1}^na_{ij}\frac{\rho_i}{2\pi n_L}\frac{\rho_j}{2\pi n_L}.$$ 
\begin{figure}[h]
\centering
\includegraphics[width=0.4\textwidth]{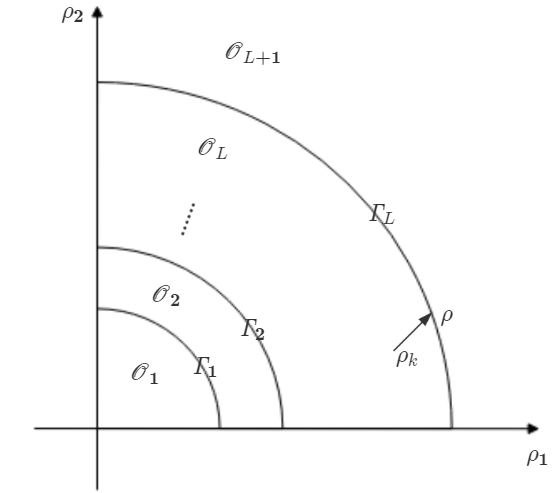}
\caption{Example of profile of bubbling solutions}
\end{figure}
We use $\mathcal{O}_{L}$ to denote the region between $\Gamma_{L-1}$ and $\Gamma_{L}$, the main purpose of this article is to study the profile of bubbling solutions $u^k$ when $\rho^k=(\rho_1^k,...,\rho_n^k)\to \rho=(\rho_1,..,\rho_n)\in \Gamma_{L}$. 
Here $\rho\in \Gamma_{L}$ ($L\ge 1$) is a limit point,
$\rho^k=(\rho_1^k,..,\rho_n^k)$ is a sequence of parameters corresponding to bubbling solutions $u^k=(u^k_1,...,u^k_n)$, which satisfies
\begin{equation}\label{b-u-sol}
\Delta_g u_i^k+\sum_j a_{ij}\rho_j^k (\frac{\mathfrak{h}_j^ke^{u_j^k}}{\int_M \mathfrak{h}_j^ke^{u_j^k}}-1)=0
\end{equation}
where  around each $p_t$ ($t=1,...,L$) $\mathfrak{h}_t^k=|x-p_t|^{2\gamma_t}h_t^k(x)$
for some smooth function $h_t^k(x)$ with uniform $C^3$ bound: There exists $c>0$ independent of $k$ such that 
\begin{equation}\label{ht-bound}
\frac 1c\le h_t^k(x)\le c, \quad \|h_t^k(x)\|_{C^3}\le c.
\end{equation}

The aim of this article is to provide a complete and precise blowup analysis for $\rho^k\to \rho$ as $k\to \infty$. The information we shall provide includes comparison of bubbling heights, difference between energy integration around blowup points, leading term of $\rho^k-\rho$, location of regular blowup points and new relation on coefficient functions around different blowup points. First we observe that 
the normal vector at $\rho$ is proportional to
$$( \sum_{j=1}^na_{1j}\frac{\rho_j}{2\pi n_L}-2,\cdots, \sum_{j=1}^na_{nj}\frac{\rho_j}{2\pi n_L}-2). $$
Based on our previous work \cite{gu-zhang} all these components are positive. For convenience we set 
\begin{equation}\label{def-m}
\mathfrak{m}_i=\sum_{j=1}^na_{ij}\frac{\rho_j}{2\pi n_L},\quad \mathfrak{m}=\min_i\{\mathfrak{m}_i\},
\end{equation}
and we assume a controlled behavior of $\rho^k\to \rho$:
\begin{equation}\label{rho-com}
\frac{\rho_i^k-\rho_i}{\rho_j^k-\rho_j}\sim 1\quad \forall  \quad i\neq j\in I:=\{1,...,n\}.
\end{equation}
Note that $A_k\sim B_k$ means $CA_k\le B_k\le C_1A_k$ for some $C,C_1>0$ independent of $k$. We say $\rho^k$ tends to $\Gamma_L$ \emph{non-tangentially} if (\ref{rho-com}) holds.

When $\rho^k\to \rho\in \Gamma_L$, if blowup occurs, let $p_1$,...,$p_{L_2}$ be blowup points, since singular sources may also be blowup points, we use $p_1.,,,p_{L_1}$ to denote regular blowup points and $p_{L_1+1}...,p_{L_2}$ to denote singular blowup points.  It is proved (see \cite{lin-zhang-jfa,gu-zhang}) that they satisfy 
\begin{equation}\label{reg-p}
\sum_{i=1}^n \bigg (\nabla (\log h_i)(p_t)+2\pi\mathrm{m}\sum_{s=1}^{L_2}\mu_s\nabla_1G^*(p_t,p_s)\bigg )\rho_i=0,\quad t=1,..,L_1.
\end{equation} 
where $\mu_s=1+\gamma_s$ (so $\mu_s=1$ if $p_s$ is a regular point), 
\begin{equation}\label{gstar}
G^*(p_t^k,p_s^k)=\left\{\begin{array}{ll}
\gamma(p_t^k,p_t^k),\quad s=t, \\
G(p_t^k,p_s^k),\quad s\neq t.
\end{array}
\right. \qquad s,t =1,...,N.
\end{equation}
Thus $(p_1,...,p_{L_1})$ is a critical point of the following function: 
\begin{align}\label{f-mul}
f(x_1,...x_{L_1})=\sum_{t=1}^{L_1}(\sum_{i=1}^n\rho_i\log h_i(x_t))+(\sum_{i=1}^n\rho_i)\pi \mathfrak{m}\sum_{t,s=1}^{L_1}G^*(x_t,x_s)\\
+(\sum_{i=1}^n\rho_i)2\pi\mathfrak{m}\sum_{t=1}^{L_1}\sum_{s=L_1+1}^{L_2}G^*(x_t,p_s).\nonumber
\end{align}
Thus if $det(D^2f)$ is not zero at each critical point, the critical points are disjoint.
  As a consequence of the blowup analysis we shall carry out in this work, we present our main theorem as follows:

\begin{thm}\label{main-thm-cur}
Let $p_1,...p_{L_1},p_{L_1+1},..,p_{L_2}$ satisfy (\ref{reg-p}) and $D^2f$ be non-degenerate at critical points of $f$. Then there exists constants $c_t^*>0$ such that if 
\begin{equation}\label{curv-c}
\Delta (\log h_i^k)(p_t)-2K(p_t^k)+2\pi n_L\ge c_t^* \quad \mbox{ for all $p_t$}, \quad t=1,...,L_1,
\end{equation}
  $u^k$ is uniformly bounded as $\rho^k$ tends to $\Gamma_L$ from above non-tangentially, where $u^k=(u_1^k,...,u_n^k)$ is a sequence of solutions to (\ref{b-u-sol}), $c_t^*$ depends only on $\mathfrak{m}$, the distance from $p_t$ to the nearest other blowup point and the injectivity radius at $p_t$, 
\end{thm}

Theorem \ref{main-thm-cur} follows naturally from the complete bubble interaction results in the next section. If the system is reduced to a single equation, a rather detailed blowup analysis has been done by many people \cite{barto2,barto3,bart-taran-jde,bart-taran-jde-2,ChenLin1,chenlin2,chen-lin-dcds,chen-lin-cpam-15,huang-zhang,kuo-lin-jdg,li-cmp,wei-zhang-adv,wei-zhang-plms, wei-zhang-jems,wei-wu-zhang,zhangcmp,zhangccm}. For a single Liouville equation, the total energy of all solutions is just one number. However, for Liouville systems defined on $\mathbb R^2$, it is established in \cite{CSW, linzhang1} that total energy of all components form a $n-1$ dimensional hyper-surface similar to $\Gamma_1$. This continuum of energy brings great difficulty to blowup analysis. A simplest description is this: if $\rho^k\to \Gamma_2$ and there are two bubbling disks. Around each bubbling disk, the local energy tends to $\Gamma_1$, but in order to know exactly how $\rho^k\to \Gamma_2$ one needs precise information about how the local energy is tending to $\Gamma_1$ and the comparison of bubbling profiles around different blowup points. This seems to be a unique feature of Liouville systems. In 2013  Lin and the second author published an article \cite{lin-zhang-jfa} that describes $\rho^k\to \Gamma_1$, which avoided this major difficulty. So the main contribution in this article is to completely present the details of bubble interactions. To put our key idea in short, we show that first, bubbling solutions are sufficiently close to global solutions in a neighborhood of the blowup points.  Since global solutions are radial, their behavior is determined by their initial conditions. In \cite{linzhang1} Lin and Zhang established an important one-to-one correspondence between the initial condition and the total integration of Liouville system. It is based on this principle we are able to obtain the precise comparison among profiles at different blowup points.  The reason in this article we require all the singular sources to have negative strength is because we have a classification theorem for such Liouville systems \cite{lin-zhang-dcds}.

The organization of this article is as follows: In section two we state our results on bubble interactions, which lead to the proof of the main theorem at the end of this section.  In section three we write the equation around a local blowup point, which happens to be a singular source. Then in section four we prove the main estimates that the profiles of bubbling solutions are extremely close around different blowup points. As a result of the main estimates of section four we prove Theorem \ref{heightthm},Theorem \ref{relation-h} and Theorem \ref{localmass}. Section five is dedicated to proving the leading terms in Theorem \ref{gammaNpneq} and Theorem \ref{main-3}. All these results depend on a local estimate proved in section six.

\section{Results on bubble interactions}

 It is established in \cite{linzhang2} that when $\rho^k=(\rho_1^k,..,\rho_n^k)$ tends to $\rho\in\Gamma_{L}$ non-tangentially, suppose there are $L_2$ blowup solutions $p_1,...,p_N$, some of them may be singular sources, some may be regular points. It is proved in \cite{gu-zhang} that 
\begin{equation}\label{relation-1}
n_L=\sum_{s=1}^{L_2} \mu_s.
\end{equation}

 In \cite{lin-zhang-jfa} Lin-Zhang derived the leading term when there is no singular source and $\rho^k$ tends to the \emph{first} critical hyper-surface $\rho^k\to \Gamma_1$. The main reason they can only discuss $\rho\to \Gamma_1$ is because there is only one blowup point in this case. Whether or not the same level of error estimate can still be obtained for multiple-bubble-situations ($\rho\to \Gamma_N$) was a major obstacle for Liouville systems in general. Recently Huang-Zhang \cite{huang-zhang-multiple} completely solved this problem by extending Lin-Zhang's result to $\rho^k\to \Gamma_N$ for any $N$. But Huang-Zhang's article does not include the singular cases. It is our  goal to further extend all the theorems in \cite{huang-zhang-multiple} to singular Liouville systems. 

One point on $\Gamma_L$ plays a particular role: Let $Q=(Q_1,...Q_n)$ be defined by $\sum_ja_{ij}Q_j=8\pi n_L$ for each $i$. Our approximating theorems depend on whether $\rho=Q$ or not. 

For simplicity we set 
\begin{equation}\label{bar-uik}
\hat u_i^k(x)=u_i^k(x)-\log \int_{M}\mathfrak{h}_i^ke^{u_i^k}dV_g,
\end{equation}
thus we have 
\[\int_M \mathfrak{h}_i^k\hat u_i^kdV_g=1,\]
and 
\begin{equation}\label{bar-u-ik}
\Delta_g \hat u_i^k(x)+\sum_{j=1}^na_{ij}\rho_j^k (\mathfrak{h}_j^ke^{\hat u_j^k}-1)=0,\quad i\in I. 
\end{equation}

Let 
\begin{equation}\label{Mk}
M_{k,t}=\max_{i\in I}\max_{x\in B(p_t,\delta_0)} \frac{\hat u_i^k(x)}{\mu_t},
\end{equation}
be the maximum over $B(p_t^k,\delta_0)$ divided by $\mu_l$. 
Here $\delta_0$ is small enough so that $B(p_{t_1},\delta_0)\cap B(p_{t_2},\delta_0)=\emptyset $ for all $t_1\neq t_2$, $\gamma_t=0$ ($\mu_t=1$) if $p_t$ is a regular point. Let
\begin{equation}\label{largest-mk}
M_k=\max_tM_{k,t},\quad \epsilon_k=e^{-\frac 12 M_k}.
\end{equation} One natural question is 
\emph{What is the relation between $M_t^k$?}. The answer of this question is presented in this theorem:
\begin{thm} \label{heightthm} Let $u^k\in \,\, \hdot^{1,n}(M)$ be a sequence of blowup solutions of (\ref{b-u-sol}). Suppose $(H1)\&(H2)$ holds for $A$, \eqref{rho-com} holds for $\rho^k$ and $\mathfrak{h}_t^k$ is described by (\ref{ht-bound}).  Then
\begin{equation}
|\mu_sM_{k,s}-\mu_t M_{k,t}|=O(1),\quad\text{ for }\,\, s,t\in\{1,\cdots, N\} .
\end{equation}
\end{thm}
\begin{rem} 
Theorem \ref{heightthm} implies that $\mu_tM_{k,t}-M_k=O(1)$ for all $t$. 
\end{rem}

It is proved in \cite{gu-zhang} that $\displaystyle{\int_{M\setminus \cup_t B(p_t,\delta_0)} \mathfrak{h}_i^ke^{\hat u_i^k}=o(1)}$. Moreover for any two blowup points $p,q$,
$$\frac{\lim_{k\to \infty}\int_{B(p,\delta)}\mathfrak{h}_i^ke^{\hat u_i^k}}{\mu_p}= \frac{\lim_{k\to \infty}\int_{B(q,\delta)}\mathfrak{h}_i^ke^{\hat u_i^k}}{\mu_q},$$
where $\delta>0$ is small to make bubbling disks mutually disjoint.
The second main question is
\emph{ What is the exact difference between them?} 

Our next result is 
\begin{thm}\label{localmass} Under the same assumptions in Theorem \ref{heightthm}.
\begin{itemize}
\item[(1)] Let
$\rho^k\to \rho\in \Gamma_L$ from $\mathcal{O}_{L}$ or $\mathcal{O}_{L+1}$ as in (\ref{rho-com}). If $\rho\neq Q$ or $\rho=Q$ but there is no regular blowup point, 
then for any $s\neq t$, 
\begin{equation}\label{each-close}
|\frac{\int_{B(p,\delta)}\mathfrak{h}_i^ke^{\hat u_i^k}}{\mu_p }-\frac{\int_{B(q,\delta)}\mathfrak{h}_i^ke^{\hat u_i^k}}{\mu_q}|=O(\epsilon_{k}^{\mathfrak{m}-2}),
\end{equation}
for $i\in I$ and $\epsilon_k=e^{-M_k/2}$, $\mathfrak{m}$ is defined in (\ref{def-m}). 
\item[(2)] Let
$\rho^k\to  Q$ from $\mathcal{O}_{L}$ or $\mathcal{O}_{L+1}$, if there is at least one regular blowup point, then \begin{equation}\label{close-m-2}
|\frac{\int_{B(p,\delta)}\mathfrak{h}_i^ke^{\hat u_i^k}}{\mu_p }-\frac{\int_{B(q,\delta)}\mathfrak{h}_i^ke^{\hat u_i^k}}{\mu_q}|=O(\epsilon_k^2\log \frac{1}{\epsilon_k}). 
\end{equation} 
\end{itemize}
\end{thm}
\begin{rem} Theorem \ref{localmass} represents the major difficulty in blowup analysis for Liouville systems. Since the mass of a system satisfies a $n-1$ dimensional hypersurface, it is particularly difficult to obtain a precise estimate on the profile of bubbling solutions around different blowup point. Theorem \ref{localmass} is a complete resolution of this long standing difficulty. \end{rem}

Then we 
define $L_2$ open sets $\Omega_{t,\delta_1}$ such that they are mutually disjoint, each of them contains a bubbling disk and their union is $M$:
\begin{equation}\label{set-decom}
B(p_t^k,\delta_1)\subset \Omega_{t,\delta_1}, \quad \cup_{t=1}^{L_2} \overline{\Omega_{t,\delta_1}}=M,
\quad \Omega_{t,\delta_1}\cap \Omega_{s,\delta_1}=\emptyset, \,\, \forall t\neq s.
\end{equation}

Next we set
\begin{align}\label{A-delta-0}
&A_{i,t,\delta}\\
=&\frac{\delta^{\mu_t(2-\mathfrak{m})}}{\mu_t}-\frac {\mathfrak{m}-2}{2\pi}\int_{\hat{\Omega}_{t,\delta}}|x-p_t|^{2\gamma_t}\frac{h_i^k(x)}{h_i^k(p_t^k)}
e^{2\pi \mathfrak{m}\sum_{l=1}^{L_2}\mu_l(G(x,p_l^k)-G^*(p_t^k,p_l^k))}dV_g \nonumber
\end{align}
where $\hat \Omega_{t,\delta}=\Omega_{t,\delta_1}\setminus B(p_t^k,\delta)$.
Then for fixed $k$ $\lim_{\delta\to 0}A_{i,\delta} $ exists
because the leading term from the  Green's function is $-\frac 1{2\pi}\log |x-p_t^k|$. In the next two theorems we identify the leading terms.

\begin{thm} \label{gammaNpneq} Suppose the same assumptions in Theorem \ref{heightthm} hold and let
$\rho^k\to \rho\in \Gamma_L$ from $\mathcal{O}_L$ or $\mathcal{O}_{L+1}$ non-tangentially. If either $\rho\neq Q$ or $\rho=Q$ but there is no regular blowup point, then
\begin{equation}\label{11july13e9}
\Lambda_L(\rho^k)
=(D+o(1))\frac{\epsilon_k^{\mathfrak{m}-2}}{n_L}.
\end{equation}
where $D$ is defined as 
$$D= \sum_{i\in I_1}\sum_{t=1}^{L_2} B_{it}\lim_{\delta\to 0}A_{i,t,\delta},$$
for $A_{i,t,\delta}$ in  (\ref{A-delta-0}) and $I_1:=\{i\in I;\quad \mathfrak{m}_i=\mathfrak{m}\}$,  $o(1)\to 0$ as $\delta\to 0$,
$\hat{\Omega}_{t,\delta_0}= \Omega _{t,\delta}\setminus B(p_t^k,\delta)$,
 and 
$$B_{it}=\frac{e^{2\pi \mathfrak{m}\sum_{l=1}^{L_2}\mu_l G^*(p_t^k,p_l^k)}}{e^{2\pi \mathfrak{m}\sum_{l=1}^{L_2}\mu_l G^*(p_1^k,p_l^k)}}\frac{h_i^k(p_t^k)}{h_i^k(p_1^k)}
e^{D_i-\alpha_i}$$ 
and $D_i$, $\alpha_i$ are constants determined by the limit of $u^k$ after scaling ( see (\ref{di-alpha})). 
\end{thm}
It can be observed that if $\rho^k$ tends to $\Gamma_L$ from $\mathcal{O}_L$ (that is, tending to $\Gamma_L$ from inside), $\Lambda_{L}(\rho^k)>0$ for $\rho^k\in \mathcal{O}_{L}$ and $\Lambda_{L}(\rho^k)<0$ for $\rho^k\in \mathcal{O}_{L+1}$. Thus, if $D\not=0$, the blowup solutions with $\rho^k$   occur   as $\rho^k\to \rho\in\Gamma_{L}$ only from one side of $\Gamma_{L}$.  Furthermore, it  yields a uniform bound of solutions as $\rho^k$ converges
to $\Gamma_{L}(\rho\not= Q)$ from $\mathcal{O}_{L}$    provided that $D<0$.   The study of the sign of $D$ is another interesting fundamental question.  Besides consideration on the compactness of solutions, the sign of $D$ is important for constructing blowup solutions with critical parameters and the uniqueness of bubbling solutions. These projects will be carried out in a different work.

The next result is concerned with the leading term of
$\Lambda_{L}(\rho^k)$ when $\rho^k\to Q$:
\begin{thm}\label{m-eq-4} Under the same assumptions in Theorem \ref{heightthm}. If   $\rho^k\to Q$ from $\mathcal{O}_{L}$ or $\mathcal{O}_{L+1}$ non-tangentially, and there is at least one regular blowup point, then
\begin{equation}\label{11july13e10}
\Lambda_{L}(\rho^k)
=-4\sum_{i=1}^n\sum_{t=1}^{L_1}b_{it}^k\epsilon_k^2\log \epsilon_k^{-1}+O(\epsilon_k^2),
\end{equation}
where 
\begin{align*}
b_{it}^k&=e^{D_i-\alpha_i}\bigg (\frac 14 \Delta (\log h_i^k)(p_t^k)-\frac{K(p_t^k)}2+2\pi N_L  \\
&+\frac 14|\nabla (\log h_i^k)(p_t^k)+8\pi\sum_{l=1}^{L_2}\nabla_1G^*(p_t^k,p_l^k)|^2
\bigg ) \nonumber
\end{align*}
and $K$ is the Gaussian curvature.
\end{thm}

\begin{rem} The leading term in Theorem \ref{gammaNpneq} is involved with integration on the whole $M$, but the leading term in Theorem \ref{m-eq-4} only depends on the information at regular blowup points. 
\end{rem}

The fifth result is about the locations of the regular blowup points:

\begin{thm}\label{main-3}
If $\rho^k\to \rho\neq Q$ as in (\ref{rho-com}), then for each regular blowup point $p_t$ ($t=1,...,L_1$)
\begin{equation}\label{11july13e7}
\sum_{i=1}^n \bigg (\nabla (\log h_i^k)(p_t^k)+2\pi \mathfrak{m}\sum_{s=1}^N \mu_s\nabla_1 G^*(p_t^k,p_s^k)\bigg )\rho_i=O(\epsilon_k^{\mathfrak{m}-2})
\end{equation}
where $\nabla_1$ means the differentiation with respect to the first component.
If $\rho^k\to Q$ as in (\ref{rho-com}),
\begin{equation}\label{11july13e8}
\sum_{i=1}^n \bigg (\nabla (\log h_i^k)(p_t^k)+8\pi \sum_{s=1}^{L_2} \mu_s\nabla_1 G^*(p_t^k,p_s^k)\bigg )q_i=O(\epsilon_k^{2}\log \epsilon_k^{-1}), 
\end{equation} for each regular blowup point $p_t$. 
\end{thm}

The sixth main result is a surprising restriction on coefficient function $h_i^k$.
\begin{thm}\label{relation-h}
Let 
$$H_{it}=\frac{2\pi\mathfrak{m}_i}{\mathfrak{m}_i-2}(\sum_{l=1}^{L_2}\mu_l G^*(p_t^k,p_l^k))+
\frac{1}{\mathfrak{m}_i-2}\log \frac{h_i^k(p_t^k)}{\mu_t^{\mathfrak{m}_i}}$$
where $\mathfrak{m}_i=\sum_{j=1}^na_{ij}\frac{\rho_j}{2\pi  n_L}$. Then 
$$H_{it}-H_{is}=H_{jt}-H_{js}+E, \quad \forall i,j=1,...,n, \quad t,s=1,...,L_2. $$
\end{thm}

Obviously Theorem \ref{relation-h} is not seen in
the case of $L_2=1$ and reveals new relations on coefficient functions $h_i^k$. In other words, if one constructs bubbling solutions, the $h_i^k$s need to satisfy the statement of Theorem \ref{relation-h}, in addition to other key information such as precise information about bubbling interactions, exact location of blowup points, accurate vanishing rate of coefficient functions and specific leading terms in asymptotic expansions. All these have been covered in the main results of this article. The construction of bubbling solutions will be carried out in other works in the future.

\subsection{Proof of Theorem \ref{main-thm-cur}}Finally in this section we provide the proof of Theorem \ref{main-thm-cur}.
We consider two cases. In the first when the the limit point $\rho\neq Q$, the blowup analysis gives $\mathfrak{m}<4$. The second case if when $\rho=Q$ and $\mathfrak{m}=4$ in this case. 

When $\rho\neq Q$ we prove that $A_{i,t,\delta}\le 0$ for each $t$ under the curvature assumption. From the definition of $A_{i,t,\delta}$ we have 
\[A_{i,t,\delta}\le \frac{\delta^{\mu_t(2-\mathfrak{m})}}{\mu_t}-\frac{\mathfrak{m}-2}{2\pi}\int_{\mathfrak{B}_{\epsilon_0}}|x-p_t|^{2\gamma_t}\frac{h_i(x)}{h_i^k(p_t)}e^{2\pi\mathfrak{m}\sum_{l=1}^{L_2}\mu_l(G(x,p_l^k)-G^*(p_t^k,p_l^k))}dV_g\]
where $\mathfrak{B}_{\epsilon_0}=B(p_t,\epsilon_0)\setminus B(p_t,\delta)$. When we evaluate the integration on $\mathfrak{B}_{\epsilon_0}$ we use $dV_g=e^{\phi}dx$ to obtain
\begin{align*}
&-\frac{\mathfrak{m}-2}{2\pi}\int_{\mathfrak{B}_{\epsilon_0}}|x-p_t|^{2\gamma_t}\frac{h_i(x)}{h_i^k(p_t)}e^{2\pi\mathfrak{m}\sum_{l=1}^{L_2}\mu_l(G(x,p_l^k)-G^*(p_t^k,p_l^k))}dV_g\\
&\le \frac{\epsilon_0^{\mu_t(2-\mathfrak{m})}}{\mu_t}-\frac{\delta^{\mu_t(2-\mathfrak{m})}}{\mu_t}\\
&-\frac{(\mathfrak{m}-2)(\epsilon_0^{2-(\mathfrak{m}-2)\mu_t}-\delta^{2-(\mathfrak{m}-2)\mu_t})}{4(2-(\mathfrak{m}-2)\mu_t)}(\Delta (\log h_i)(p_t)-2K(p_t)+2\pi n_l)
\end{align*}
Note that the reason we have an inequality is because we drop the square of the first derivatives in the exponential function, which is positive. Thus, as long as
\[\Delta (\log h_i)(p_t)-2K(p_t)+2\pi n_l>\frac{4(2-(\mathfrak{m}-2)\mu_t)}{(\mathfrak{m}-2)\mu_t} \epsilon_0^{-2},\quad \forall t,\]
 $\Lambda_{L}(\rho^k)$ is negative, so if $\rho^k$ is tending to $\rho$ from above non-tangentially, the blow-up does not happen. The existence of $\epsilon_0$ depends on the injectivity radius, the distance from $p_t^k$ to other blowup point and singularities. The argument for $\rho^k\to Q$ is similar.  Theorem \ref{main-thm-cur} is established. $\Box$

\section{The profile of bubbling solutions around a singular source}
In this section we write the equation (\ref{b-u-sol}) around a singular source $p$, which is also a blowup point. Suppose that the strength of the singularity is $4\pi \gamma_p$ with $\gamma_p\in (-1,0)$, we derive an approximation theorem for blow-up solutions $u^k$ around $p$. 
To state more precise approximation results we write the equation in local coordinate around $p$. In this coordinate,  $ds^2$ has the form

$$e^{\phi(y_{p)}}[(dy^1)^2+(dy^2)^2],$$
where
$$\nabla\phi(0)=0,\phi(0)=0.$$
Also near $p$ we have
$$\Delta_{y_{p}}\phi=-2Ke^\phi,\quad {\rm where}\hspace{0.1cm} K {\rm \hspace{0.1cm}is\hspace{0.1cm} the \hspace{0.1cm}Gauss\hspace{0.1cm} curvature}. $$
 Here we invoke a result proved in our previous work \cite{gu-zhang} since $\gamma_p$ is not a positive integer, the spherical Harnack inequality holds around $p$ and $p$ is the only blowup point in $B(p,\delta)$. Here we recall that in a neighborhood of $p$, $\mathfrak{h}_i^k(x)=|x|^{2\gamma_p}h_i^k(x)$.

In this local coordinates, (\ref{b-u-sol}) is of the form

\begin{equation}\label{4}
-\Delta u_i^k=\sum_{j=1}^n a_{ij}\rho_j^ke^{\phi}(\frac{|x|^{2\gamma_p}h_j^{k}e^{u_j^k}}{\int_M\mathfrak{h}^k_je^{u_j^k}}-1),\hspace{0.2cm} B(0,\delta).
\end{equation}

Going back to (\ref{4}), we let $f_i^k$ be defined as

$$ \Delta f_i^k=\sum_{j=1}^n a_{ij}\rho_j^k e^\phi,\quad {\rm in}\hspace{0.1cm}B(0,\delta),\quad \mbox{and}\quad 
f_i^k(0)=|\nabla f_i^k(0)|=0. $$
Then we have
$$ \Delta (\hat u_i^k-f_i^k)+\sum_{j=1}^n a_{ij}\rho_j^k |x|^{2\gamma_p} h_j^{k}e^{\hat u_j^k-f_j^k}e^{f_j^k}e^\phi=0,\quad {\rm in}\hspace{0.1cm}B(0,\delta),$$
which can be further written as 
\begin{equation}\label{6}
\Delta \tilde{u}_i^k+\sum_{j=1}^n a_{ij}|x|^{2\gamma_p}\tilde{h}_j^{k}e^{\tilde{u}_j^k}=0,\hspace{0.2cm} B(0,\delta).
\end{equation}
if we set
$$\tilde{h}_i^k(x)=\frac{h_i^{k}(x)}{h_i^{k}(p)}e^{\phi+f_i^k},
\quad \mbox{which implies} \quad \tilde{h}_i^k(0)=1, $$
and 
\begin{equation}\label{5}
\tilde{u}_i^k=\hat u_i^k+\log(\rho_i^k h_i^{k}(p))-f_i^k,
\end{equation}

Now we introduce $\phi_i^k$ to be a harmonic function defined by the oscillation of $\tilde{u}_i^k$ on $B(p,\delta)$:

\begin{equation}\label{7}
\left\{\begin{array}{ll}
  -\Delta \phi_i^k=0,\quad {\rm in}\hspace{0.1cm}B(0,\delta), \\
  \phi_i^k=\tilde{u}_i^k-\frac{1}{2\pi\delta}\int_{\partial B(0,\delta)}\tilde{u}_i^k, \quad {\rm on}\hspace{0.1cm}\partial B(0,\delta)
  \end{array}\right.
\end{equation}

Obviously, $\phi_i^k(0)=0$ by the mean value theorem and $\phi_i^k$ is uniformly bounded on $B(0,\delta/2)$ because $u_i^k$ has finite oscillation away from blowup points.
Now we set
\begin{equation}\label{local-M}M_{k,p}=\frac{\max_{i\in I} \tilde u_i^k(p)}{\mu_p}, \quad \varepsilon_{k,p}=e^{-\frac{M_{k,p}}{2}}, \quad \mu_p=1+\gamma_p
\end{equation} and 
$$
\sigma_{ip}^k=\frac 1{2\pi}\int_{B(p,\delta)}|x-p_k|^{2\gamma_p}\tilde h_i^k e^{\tilde u_i^k},\quad m_{ip}^k=\sum_{j=1}^n a_{ij}\sigma_{ip}^k,\quad m_p^k=\min_{i\in I}m_{ip}^k. $$
It is proved in \cite{gu-zhang} that 
$$v_i^k(y)=\tilde u_i^k(p+\epsilon_{k,p}y)+2(1+\gamma_p)\log \epsilon_{k,p}, \quad |y|<\delta \epsilon_{k,p}^{-1} $$
converges uniformly to $U=(U_1,....,U_n)$ over any fixed compact subset of $\mathbb R^2$ and $v$ satisfies 
\begin{equation}\label{v-global}
\Delta U_i+\sum_j a_{ij}|y|^{2\gamma_p}e^{U_j}=0,\quad \mbox{in}\quad \mathbb R^2, \quad i=1,..,n. 
\end{equation}
In other words, after the scaling, no component is lost in the limit. We have a fully bubbling sequence. It is also established in \cite{gu-zhang} that $p$ is the only blowup point in a neighborhood of $p$ and $u^k$ satisfies spherical Harnack inequality around $p$. 

Let $U^k=(U_1^k,\cdots,U_n^k)$ be the radial solutions of
\begin{equation}\label{8}
\left\{\begin{array}{ll}
  -\Delta U_i^k=\sum_{j=i}^n a_{ij}|x|^{2\gamma_p} e^{U_j^k},\quad {\rm in}\hspace{0.1cm}R^2, \\
  U_i^k(0)=v_i^k(0), \quad i\in I.
  \end{array}\right.
\end{equation}
This family of global solutions $U^k$ will be used as the first term in the approximation of $v^k$. 

Based on Theorem \ref{2.1} we have the following estimate of $v_i^k-U_i^k$: For multi-index $\alpha$ ($|\alpha|=0$ or $1$)
\begin{align*}
    &|D^{\alpha}(v_i^k(y)-U_i^k(y)-\phi_i^k(\epsilon_ky)-\Phi_i^k(y)|\\
    \le & C\epsilon_{k,p}^2(1+|y|)^{2+2\mu_p-m_p-|\alpha |+\epsilon} \quad |y|\le \tau \epsilon_{k,p}^{-1}
\end{align*}
where $\phi_i^k$ is defined in (\ref{7}), $C>0$ is also independent of $k$, $m_p=\lim_{k\to \infty} m_p^k$,
$$\Phi_i^k(y)=\epsilon_{k,p}(G_{1,i}^k(r)\cos \theta+G_{2,i}^k(r)\sin \theta), $$
with 
$$|G_{t,i}^k(r)|\le Cr(1+r)^{2\mu_p-m_p+\epsilon} $$
where $\epsilon>0$ is a small positive constant. 

\section{The comparison of bubbling profiles}

For simplicity 
we first assume there are only two blowup points $p$ and $q$. The conclusion for the more general situation will be stated at the end of this section. If $p$ is a singular source, we set 
$M_{k,p}$ be defined as in (\ref{local-M}).
If $p$ is a regular point, we set 
$$M_{k,p}=\max_{x\in B(p,\delta)}\max_i \tilde u_i^k(x),$$
for some small $\delta>0$. $M_{k,q}$ is understood in the same fashion. Here we require that at least one of $p,q$ is a singular source because the comparison of bubbling profiles for regular blowup points have been done in \cite{huang-zhang-multiple}. 
Now we use Green's representation to describe the neighborhood of $p$. By (\ref{bar-u-ik})  $\hat u_i^k$ is
\begin{equation*}
\begin{aligned}
\hat u_i^k(x)&=\Bar{u}_i^k+\int_M G(x,\eta)\sum_{j=1}^n a_{ij} \rho_j^k \mathfrak{h}_j^k e^{\hat u_j^k}dV_g\\
&=\Bar{u}_i^k+\left(\int_{B(p,\delta)}+\int_{B(q_k,\delta)}+\int_{M\setminus B(p,\delta)}\right)G(x,\eta)\sum_{j=1}^n a_{ij} \rho_j^k \mathfrak{h}_j^k e^{\hat u_j^k}\\
&=\Bar{u}_i^k+\uppercase\expandafter{\romannumeral1}+\uppercase\expandafter{\romannumeral2}+\uppercase\expandafter{\romannumeral3}.
\end{aligned}
\end{equation*}

Here 
$$\Bar{u}_i^k=\int_M \hat u_i^k dV_g, \quad E=O(\epsilon_k^{\tau})$$ for some $\tau>0$. Note that it is already established in our previous work \cite{gu-zhang} that $(\mu_pM_{k,p})/(\mu_q M_{k,q})\to 1$ as $k\to \infty$, the error term around $p$ in local approximation and computation of Pohozaev identity is $O(\epsilon_{k,p}^{m_{k,p}-2\mu_p})$ and the error around $q$ is the similar expression with $q$ replacing $p$ and with a possible correction of a logarithmic term.
Right now we are using crude bound for the error $E$. This bound will be improved later. 

For simplicity we use
$$m_i^k=\sum_{j=1}^n a_{ij}\sigma_j^k, \quad \sigma_i^k=\frac{1}{2\pi}\int_{B(p_k,\delta)} \rho_i^k h_i^k |\eta|^{2\gamma_p}e^{\hat u_i^k}dV_g$$
and $\bar m_i^k, \Bar \sigma_i^k$ are the integration in $B(q_k,\delta)$.
Also we use $m_k$ to denote the smallest $m_i^k$ and $\bar m_k$ is the smallest $\bar m_i^k$. Using Theorem \ref{2.1} to evaluate integrals we have
\begin{equation*}
\begin{aligned}
\uppercase\expandafter{\romannumeral1}&=\int_{B_\delta} (-\frac{1}{2\pi}\log|\eta|+\nu(x,\eta))\sum_{j\in I}a_{ij}|\eta|^{2\gamma}\tilde h_j^k e^{\tilde u_j^k}d\eta\\
&=-m_i^k\log|x|+2\pi m_i^k \nu(x,p_k)+E\\
\uppercase\expandafter{\romannumeral2}&=2\pi \Bar m_i^k G(x,q_k)+E,\quad 
\uppercase\expandafter{\romannumeral3}=E
\end{aligned}
\end{equation*}

Thus in the neighborhood of $p$,
\begin{equation}\label{3.1}
\hat u_i^k(x)=\Bar{u}_i^k-m_{i}^k\log|x|+2\pi m_{i}^k\nu(x,p_k)+2\pi \bar m_i^k G(x,q_k)+E
\end{equation}
and by Lemma 2.3 of \cite{gu-zhang}, we have
\begin{equation}\label{ubar}
\bar u_i^k(x)=-\frac{m_i-2\mu_p}{2}M_{k,p}+O(1).
\end{equation}

\noindent{\bf Proof of Theorem \ref{heightthm}}. 

From (\ref{ubar}), we obtain
\begin{equation}\label{l1}
-\frac{m_i-2\mu_p}{2}M_{k,p}=-\frac{\bar m_i-2\mu_q}{2}M_{k,q}+O(1).
\end{equation}

Let
$$\lambda_k=\frac{\mu_p M_{k,p}}{\mu_q M_{k,q}},\quad \delta_i^k=\frac{O(1)}{M_{k,q}},$$
we can write (\ref{3.1}) as 
\begin{equation}\label{l2}
-\frac{m_i-2\mu_p}{2\mu_p}\lambda_k+\delta_i^k=-\frac{\bar m_i-2\mu_q}{2\mu_q}.
\end{equation}
Since the Pohozaev identity gives
$$\sum_{i,j}a_{ij}\sigma_i^k\sigma_j^k=4\mu_p\sum_i \sigma_i^k+E,$$
$m^k=(m_1^k,\cdots,m_n^k)$ satisfies
\begin{equation}\label{l3}
\sum_{i,j}a^{ij}\bigg( \frac{m_i^k-2\mu_p}{2\mu_p}\bigg)\bigg( \frac{m_j^k-2\mu_p}{2\mu_p}\bigg)=\sum_{i,j}a^{ij}+E.
\end{equation}
and (\ref{l3}) also holds for $\bar m^k=(\bar m_1^k,\cdots,\bar m_n^k)$. Using (\ref{l2}) we have 
\begin{equation}\label{l4}
\sum_{i,j}a^{ij}\bigg( \frac{m_i^k-2\mu_p)}{2\mu_p}\lambda_k+\delta_i^k\bigg)\bigg( \frac{m_j^k-2\mu_p)}{2\mu_p}\lambda_k+\delta_j^k\bigg)=\sum_{i,j}a^{ij}+E.
\end{equation}
which can be written as
\begin{equation}\label{l5}
\begin{aligned}
&\lambda_k^2  \sum_{i,j}a^{ij}\bigg( \frac{m_i^k-2\mu_p}{2\mu_p}\bigg)\bigg( \frac{m_j^k-2\mu_p}{2\mu_p}\bigg)\\
&+2\lambda_k\sum_{i,j}a^{ij}\bigg( \frac{m_i^k-2\mu_p}{2\mu_p}\bigg)\delta_j^k
+\sum_{i,j}a^{ij}\delta_i^k\delta_j^k
=\sum_{i,j}a^{ij}+E.
\end{aligned}
\end{equation}
Let
$$B_k=\frac{\sum_{i,j}a^{ij} \bigg( \frac{m_i-2\mu_p}{2\mu_p}\bigg)\delta_j^k}{\sum_{i,j}a^{ij}},\quad C_k=\frac{\sum_{i,j}a^{ij}\delta_i^k\delta_j^k}{\sum_{i,j}a^{ij}}$$
be the coefficients of the following polynomial for $\lambda_k$: 
$$\lambda_k^2+B_k\lambda_k+C_k=1+E.$$
Our goal is to obtain an upper bound of $|\lambda_k-1|$. Here we note that since $A$ satisfies $(H1)$ and $(H2)$, $\sum_{i,j}a^{ij}\neq 0$.
It is obvious that $\lim_{k\rightarrow \infty} \lambda_k=1$. Taking advantage of $C_k=O( M_{k,q}^{-2})$ and $B_k=O( M_{k,q}^{-1})$ we have
\begin{equation}\label{l6}
\lambda_k=-\frac{B_k}{2}+\sqrt{\frac{B_k^2}{4}-(C_k-1)}=1+O(M_{k,q}^{-1}).
\end{equation}
This verifies $\mu_p M_{k,p}-\mu_q M_{k,q}=O(1)$. 
Theorem \ref{heightthm} is established. $\Box$

\begin{rem}
Since (see \cite{gu-zhang})
$$\lim_{k\to \infty}\frac{m_i^k}{\mu_p}
=\lim_{k\to\infty}\frac{\bar m_i^k}{\mu_q}$$
and $\mu_p M_{k,p}-\mu_q M_{k,q}=O(1)$, we have $M_{k,p}=O(\log 1/\epsilon_{k,p})$, 
and a more precise description of $O(\epsilon_{k,p}^{m_k-2\mu_p})$ and $O(\epsilon_{k,q}^{\bar m_k-2\mu_q})$:
$$\epsilon_{k,p}^{m_k-2\mu_p}=e^{-\frac 12\mu_pM_{k,p}(m_{k,p}/\mu_p-2)}
=O(1)e^{-\frac 12\mu_qM_{k,q}(\bar m_k/\mu_q-2)}=O(\epsilon_{k,q}^{\bar m_k-2\mu_q}). $$
Thus
$O(\varepsilon_{k,p}^{m_k-2\mu_p})=O(\varepsilon_{k,q}^{\bar m_k-2\mu_q})=E$. From now on in this section we use $E$ to denote $O(\epsilon_{k,p}^{m_k-2\mu_p})$. 
\end{rem}

\medskip

By the definition of $\tilde u_i^k$ in (\ref{5}) and Theorem \ref{2.1}, the value of $\tilde u_i^k$ away from bubbling disks is
\begin{equation}\label{3.2}
\begin{aligned}
\tilde u_i^k(x)=&\Bar{u}_i^k-m_i^k\log|x|+2\pi m_i^k\nu(x,p_k)+2\pi \Bar m_i^k G(x,q_k)\\
&+\log(\rho_i^k h_i^k(p_k))-f_i^k(x)+E
\end{aligned}
\end{equation}
On the other hand  $\tilde u_i^k$ also has the form
\begin{equation}\label{3.3}
\tilde u_i^k(x)-\phi_i^k(x)=V_i^k(x)+E
\end{equation}
where $V_i^k$ is the global radial solution that takes the initial value of $\tilde u^k$: $V_i^k(0)=\tilde u_i^k(p_k)$. 
From this expression we see that $\phi_i^k$ denotes the ``non-radial part" in the expansion of $\tilde u_i^k$. Combining (\ref{3.2}), (\ref{3.3}) and using the fact that $\phi_i^k(0)=0$, we have,
\begin{equation}\label{3.4}
\begin{aligned}
\phi_i^k(x)=&2\pi m_i^k(\nu(x,p_k)-\nu(p_k,p_k))+2\pi \Bar m_i^k (G(x,q_k)-G(p_k,q_k))\\
&-f_i^k(x)+E
\end{aligned}
\end{equation}

So in the local coordinate around $p$, $\tilde u_i^k$ can be rewritten as
\begin{equation}\label{3.5}
\begin{aligned}
\tilde u_i^k(x)&=-m_i^k\log|x|+\Bar{u}_i^k+\phi_i^k(x)\\
&+2\pi m_i^k \nu(p_k,p_k)+2\pi \Bar m_i^k G(p_k,q_k)+\log(\rho_i^k h_i^k(p_k))+E
\end{aligned}
\end{equation}
for $x\in B(0,\delta)\setminus B(0,\frac{\delta}2)$. 
Thus around $p_k$ we have
\begin{equation}\label{3.6}
\begin{aligned}
V_i^k(x)&=-m_i^k\log|x|+\Bar{u}_i^k\\
&+2\pi m_i^k \nu(p_k,p_k)+2\pi \Bar m_i^k G(p_k,q_k)+\log(\rho_i^k h_i^k(p_k))+E
\end{aligned}
\end{equation}

Similarly around $q_k$ we have
\begin{equation}\label{3.7}
\begin{aligned}
\bar V_i^k(x)&=-\Bar m_i^k\log|x|+\Bar{u}_i^k\\
&+2\pi\Bar m_i^k \nu(q_k,q_k)+2\pi  m_i^k G(q_k,p_k)+\log(\rho_i^k h_i^k(q_k))+E
\end{aligned}
\end{equation}

In order to obtain precise  estimate between $m_i^k$ and $\bar m_i^k$, we need the following asymptotic estimate of a global solution:

\begin{lem}\label{global-U-asy}
 Let $U=(U_1,\cdots,U_n)$ be the global solution of
\begin{equation*}
\Delta U_i+\sum_{j=1}^n a_{ij}|y|^{2\gamma}e^{U_j}=0, \int_{\mathbb R^2}|y|^{2\gamma} e^{U_j}<\infty
\end{equation*}
where $\gamma>-1$, $A=(a_{ij})_{n\times n}$ satisfies $(H1)$,
$U_i$ is radial and suppose 
$$\max_{i\in I}U_i(0)=0.$$ Then using $\mu=1+\gamma$,
\begin{equation}\label{15}
\begin{aligned}
U_i(r)&=-m_{iu}\log r+D_{iu}-\alpha_{iu}\\
&-\sum_{j=1}^n\frac{a_{ij}}{(m_{ju}-2\mu)^2}e^{D_{ju}-\alpha_{ju}}r^{2\mu-m_{ju}}+O(r^{2\mu-m_u-\delta}),\quad r>1,
\end{aligned}
\end{equation}
where $\delta>0$ and 
\begin{equation}\label{sigma-m-u}
\sigma_{iu}=\frac{1}{2\pi}\int_{\mathbb R^2}|y|^{2\gamma}e^{U_i},\quad m_{iu}=\sum_{j=1}^n a_{ij}\sigma_{ju}
\end{equation}
\begin{equation}\label{d-alpha}
D_{iu}=\int_0^{\infty}\log r \sum_{j=1}^n a_{ij}r^{2\gamma}e^{U_j(r)}rdr, \qquad \alpha_{iu}=-U_i(0),
\end{equation}
and $m_u=\min_i m_{iu}$,
\begin{equation}\label{16}
\sigma_{iu}=\sigma_{iR}+\frac{e^{D_{iu}-\alpha_{iu}}}{m_{iu}-2\mu}R^{2\mu-m_{iu}}+O(R^{2\mu-m_u-\delta})
\end{equation}
where
$$ \hspace{0.1cm}\sigma_{iR}=\frac{1}{2\pi}\int_{B_R}|y|^{2\gamma}e^{U_i}.$$
\end{lem}
\begin{rem} It is proved in \cite{lin-zhang-dcds} that if $\gamma\in (-1,0)$, all $U_i$s are radial functions. 
\end{rem}

\begin{rem}\label{pohozaev-R}
From (\ref{sigma-m-u}) and the Pohozaev identity for $(\sigma_{1u},....,\sigma_{nu})$ we have
$$4\sum_i\frac{\sigma_{iR}}{\mu}-\sum_{ij}a_{ij}\frac{\sigma_{iR}}{\mu}\frac{\sigma_{jR}}{\mu}=
2\sum_i
\frac{e^{D_{iu}-\alpha_{iu}}}{\mu^2}R^{\mu(2-m_{iu}/\mu)}+O(R^{2\mu-m_{iu}-\delta}). 
$$
\end{rem}
\noindent{\bf Proof of Lemma \ref{global-U-asy}:}
It is well known that 
\begin{equation*}
 U_i(x)=-\frac{1}{2\pi}\int_{\mathbb R^2} \log|x-y|(\sum_{j=1}^na_{ij}|y|^{2\gamma}e^{U_j})dy+c_i.
\end{equation*}
Thus 
\begin{equation*}
-\alpha_{iu}=U_i(0)=-\int_0^\infty \log r\sum_{j=1}^n a_{ij} r^{2\gamma} e^{U_j(r)}rdr+c_i=-D_{iu}+c_i.
\end{equation*}
Therefore $c_i=D_{iu}-\alpha_{iu}$ and 
\begin{equation}\label{U-exp-2}
\begin{aligned}
U_i(x)&=-\frac{1}{2\pi}\int_{\mathbb{R}^2}(\log|x-y|-\log|x|)\sum_{j=1}^n a_{ij}|y|^{2\gamma}e^{U_j(y)}dy\\
&+D_{iu}-\alpha_{iu}-m_{iu}\log|x|\\
&=-m_{iu}\log|x|+D_{iu}-\alpha_{iu}+O(r^{-\delta}),\quad r>1
\end{aligned}
\end{equation}
for some $\delta>0$. Here we recall that $m_{iu}$ is defined in (\ref{sigma-m-u}). This expression gives
$$e^{U_i(r)}=r^{-m_{iu}}e^{D_{iu}-\alpha_{iu}}+O(r^{-m_u-\delta}).$$
Then
\begin{equation*}
\begin{aligned}
\sigma_{iu}&=\sigma_{iR}+\frac{1}{2\pi}\int_{\mathbb R^2\setminus B_R}|y|^{2\gamma}e^{U_i}\\
&=\sigma_{iR}+\frac{e^{D_{iu}-\alpha_{iu}}}{m_{iu}-2\mu}R^{2+2\gamma-m_{iu}}+O(R^{2\mu-m_u-\delta}).
\end{aligned}
\end{equation*}

Since $U_i$ satisfies the following ordinary differential equation:
\begin{equation}\label{17}
U_i''(r)+\frac{1}{r}U_i'(r)=-\sum_{j=1}^n a_{ij}r^{2\gamma}e^{U_j},\quad 0<r<\infty.
\end{equation}
Multiplying $r$ on both sides and using the fact 
$\lim_{r\rightarrow \infty}rU_i'(r)=-m_{iu}$ we have
\begin{equation*}
\begin{aligned}
-m_{iu}-rU_i'(r)&=-\sum_{j=1}^n a_{ij} \int_r^\infty s^{2\gamma+1} e^{U_j(s)}ds\\
&=-\sum_{j=1}^\infty a_{ij} \frac{e^{D_{ju}-\alpha_{ju}}}{m_{ju}-2\mu} r^{2\mu-m_{ju}}+O(r^{2\mu-m_u-\delta}).
\end{aligned}
\end{equation*}
Then we have
$$U_i'(r)=-\frac{m_{iu}}{r}+\sum_{j=1}^n \frac{a_{ij}}{m_{ju}-2\mu}e^{D_{ju}-\alpha_{ju}}r^{1+2\gamma-m_{ju}}+O(r^{1+2\gamma-m_u-\delta}).$$
After integration we obtain (\ref{15}) from comparison with (\ref{U-exp-2}).
Lemma \ref{global-U-asy} is established. $\Box$

\medskip

\begin{prop}\label{close-m} If $\lim_{k\to \infty}\frac{m_i^k}{\mu_p}<4$ or $\lim_{k\to \infty}\frac{m_i^k}{\mu_p}=4$ and $\gamma_p,\gamma_q\neq 0$ ($\mu_p=1+\gamma_p$),
$$ \bigg|\frac{m_i^k}{\mu_p}-\frac{\Bar m_i^k}{\mu_q}\bigg|\leq C\delta_0^{2+2\mu_p-m_k}\varepsilon_k^{m_k-2\mu_p}, \quad i=1,...,n$$
where $m_k=\min_i m_{ik}$ 
\end{prop}

\begin{rem} Proposition \ref{close-m} is essential for all the main results. The rough idea of the proof is that both $m_i^k$ and $\bar m_i^k$ are very close to the energy of the approximating global solutions. The comparison of the global solutions at different blowup points is based on a crucial estimate of initial-value-dependence of Liouville systems established by Lin-Zhang \cite{lin-zhang-dcds}.
\end{rem}
\begin{rem}
It is trivial but important to observe that $\delta_0^{2+2\mu_p-m_k}\to 0$ if $\delta_0\to 0$.
\end{rem} 
\noindent{\bf Proof of Proposition \ref{close-m}.:}
Let $V^k=(V_1^k,...,V_n^k)$ be a sequence of global solutions that satisfies $V_i^k(p_k)=\tilde u_i^k(p_k)$. Since $\gamma_p\le 0$, $V^k$ is radial and is the first term in the approximation of $\tilde u^k$ around $p_k$. On the other hand 
$\bar V^k=(\bar V_1^k,...,\bar V_n^k)$ is the first term in the approximation of $\tilde u^k$ at $q_k$.  We use $m_{iv}^k$ and $\Bar m_{iv}^k$ to denote the integration of the global solution $V_i^k, \Bar V_i^k$ respectively. 

First we observe that the expansion of bubbling solutions in Theorem \ref{2.1} gives 
$$m_i^k-m_{iv}^k=E$$

If we use
$$U_i^k(y)=V_i^k(x)+2\mu_p\log \varepsilon_{k,p},$$
where
$$\varepsilon_{k,p}=e^{-\frac{M_{k,p}}{2}},\quad x=\varepsilon_{k,p} y$$
Then based on the expansion of $U_i^k$ in Lemma \ref{global-U-asy} we have
\begin{equation}\label{3.8}
V_i^k(x)=-m_{iv}^k\log|x|-\frac{m_{iv}^k-2\mu_p}{2}M_{k,p}+D_i-\alpha_i+E.
\end{equation}
where 
\begin{equation}\label{di-alpha}
D_i=\int_0^{\infty}\log r\sum_{j=1}^na_{ij}r^{1+2\gamma_p}e^{U_j^k}dr,
\quad \alpha_i=-U_i^k(0). 
\end{equation}
Similarly we have
\begin{equation}\label{3.9}
\Bar V_i^k(x)=-\Bar m_{iv}^k\log|x|-\frac{\Bar m_{iv}^k-2\mu_q}{2} M_{k,p}+\Bar D_i-\Bar \alpha_i+E.
\end{equation}
Note that it is proved in \cite{gu-zhang} that $u^k$ is a fully bubbling sequence around each blowup point, which means after scaling no component is lost in the limit system. The assumption (H2) is important for this result. Since $V^k=(V_1^k,...,V_n^k)$ and $\bar V^k=(\bar V_1^k,...,\bar V_n^k)$ satisfy different ODE systems, we need to make the following change of variable before comparing them.
Let

$$\tilde V_i^k(r)=\bar V_i^k(r^{\frac{\mu_p}{\mu_q}})+2\log\frac{\mu_p}{\mu_q}$$
Here we first make an important observation: If we use $\bar \alpha=(\bar \alpha_1,...,\bar \alpha_n)$ as the initial value of $\bar V^k$ after scaling, we see immediately that $\bar \alpha$ is also the initial condition of $\tilde V^k$ after scaling. This plays a crucial role in the proof of Proposition \ref{close-m}.

By direct computation we see that $\tilde V^k=(\tilde V_1^k,...,\tilde V_n^k)$ satisfies the same equation as $V^k$:
$$\frac{d^2}{dr^2}\tilde V_i^k(r)
+\frac 1r \frac{d}{dr}\tilde V_i^k(r)+\sum_j a_{ij}r^{2\gamma_p}e^{\tilde V_j^k(r)}=0,\quad 0<r<\infty, $$
and the asymptotic expansion of $\bar V_i^k$ in (\ref{3.7}) yields
\begin{equation}\label{tilde-v}
\tilde V_i^k(x)=-\frac{\mu_p}{\mu_q}\Bar m_{iv}^k\log|x|-\frac{\Bar m_{iv}^k-2\mu_q}{2}  M_{k,q}+\bar D_i-\bar \alpha_i+2\log\frac{\mu_p}{\mu_q}+E.
\end{equation}

We observe that the height around $p$ is $\mu_p M_{k,p}$ and the height around $q$ for $\tilde V^k$ is $\mu_q M_{k,q}+2\log \frac{\mu_p}{\mu_q}$. It is critical to re-scale $\tilde V^k$ to make them match.
Let
\begin{equation}\label{height}
2\mu_p\log\eta=\mu_p M_{k,p}-\mu_q M_{k,q}-2\log\frac{\mu_p}{\mu_q}
\end{equation}
and set 
\begin{equation}\label{hat-v-d}
\hat V_i^k(r)=\tilde V_i^k(\eta r)+2\mu_p\log \eta.
\end{equation}
Then the expansion of $\tilde V_i^k(x)$ in (\ref{tilde-v}) and ($\ref{height}$) lead to the following expansion of $\hat V_i^k$: 
\begin{equation}\label{3.10}
\begin{aligned}
\hat V_i^k(x)&=\tilde V_i^k(\eta x)+2\mu_p\log\eta\\
=&-\frac{\mu_p}{\mu_q}\Bar m_{iv}^k\log|x|-\frac{\Bar m_{iv}^k-2\mu_q}{2\mu_q} M_{k,p}\mu_p\\
&+\bar D_i-\bar \alpha_i+\frac{\Bar m_{iv}^k}{\mu_q}\log\frac{\mu_p}{\mu_q}+E.
\end{aligned}
\end{equation}

As mentioned before if we use $-\hat \alpha_i=\hat U_i(0)$ we see that $\hat \alpha_i=\bar \alpha_i$ for all $i$. This is also due to the fact that $\hat \alpha_i$ is the difference between the initial value with the largest initial value. Since from $\tilde V_i$ to $\hat V_i$ all components are added by a same number, the difference between any two components remains the same. 
Since both $V_i^k$ and $\hat V_i^k$ are radial and satisfy the same Liouville system, the dependence on initial condition gives
\begin{equation}\label{3.11}
|V_i^k(x)-\hat V_i^k(x)|\leq C\sum_{i=1}^n|\alpha_i-\bar \alpha_i|
\end{equation}

Suppose  $\alpha_1=\bar \alpha_1=0$ and the mapping from $(\alpha_2,\cdots,\alpha_n)$ to $(\sigma_2,\cdots,\sigma_n)$ is a diffeomorphism ( see \cite{lin-zhang-dcds}). Let 
$$ \hat \sigma_i^k=\frac{1}{2\pi}\int_{\mathbb R^2} |y|^{2\gamma_p}e^{\hat V_i^k}, \quad i=1,...,n . $$
Using the definition of $\hat V^k$ in (\ref{hat-v-d}) we see that 
$$\hat \sigma_i^k=\frac{\mu_p}{\mu_q}\bar \sigma_i^k. $$
Thus (\ref{3.11}) can be written as
\begin{equation}\label{3.12}
|V_i^k(x)-\hat V_i^k(x)|\le C\sum_{i}^n\bigg|\frac{\sigma_i^k}{\mu_p}-\frac{\Bar \sigma_i^k}{\mu_q}\bigg|
\end{equation}

From (\ref{3.8}) and (\ref{3.10}) we see that the difference between $V_i^k$ and $\hat V_i^k(x)$ is

\begin{equation}\label{3.13}
\begin{aligned}
&V_i^k(x)-\hat V_i^k(x) \\
=&(\frac{\mu_p}{\mu_q}\Bar m_{iv}^k-m_{iv}^k)\log|x|
+\bigg(\frac{\Bar m_{iv}^k}{2\mu_q}-\frac{m_{iv}^k}{2\mu_p}\bigg)M_{k,p}\mu_p\\
&-(\alpha_i-\bar\alpha_i)+(D_i-\bar D_i-\frac{\Bar m_{iv}^k}{\mu_q}\log\frac{\mu_p}{\mu_q})+E.
\end{aligned}
\end{equation}

Here we point out another key fact: The maximum heights of $V^k$ and $\hat V^k$ are the same. So when we scale $\hat V^k$ to $\hat U^k=(\hat U_1^k,...,\hat U_n^k)$, the scaling factor is still $\epsilon_{k,p}$:

\begin{equation*}
\begin{aligned}
\hat U_i^k(x)&=\hat V_i^k(\varepsilon_{k,p} x)+2\mu_p\log \varepsilon_{k,p}\\
&=\tilde V_i^k(\eta \varepsilon_{k,p} x)+2\mu_p\log \eta+2\mu_p\log\varepsilon_{k,p}\\
&=\bar V_i^k((\eta \varepsilon_{k,p} x)^{\frac{\mu_p}{\mu_q}})+2\log\frac{\mu_p}{\mu_q}+2\mu_p\log(\eta\varepsilon_{k,p})
\end{aligned}
\end{equation*}
Let $\hat D_i$ be defined by 
\begin{equation*}
\begin{aligned}
\hat D_i&=\int_0^\infty (\log r) \, r^{2\gamma_p+1}\sum_j a_{ij}e^{\hat U_j^k(r)}dr\\
&=\int_0^\infty (\log r) r^{2\gamma_p+1}\sum_j a_{ij}e^{\bar V_j^k((\eta\varepsilon_{k,p} r)^{\frac{\mu_p}{\mu_q}})}(\frac{\mu_p}{\mu_q})^2 \eta^{2\mu_p}\varepsilon_{k,p}^{2\mu_p} dr
\end{aligned}
\end{equation*}
We choose $s$ to make 
$$\eta \epsilon_{k,p} r= \epsilon_{k,q} s $$
where $\epsilon_{k,p}=e^{-\frac 12 M_{k,p}}$, $ \epsilon_{k,q}=e^{-\frac 12  M_{k,q}}$. From (\ref{height}) and 
direct computation
$$s=(\frac{\mu_q}{\mu_p})^{\frac{1}{\mu_q}}r^{\frac{\mu_p}{\mu_q}}$$
With this change of variable we evaluate $\hat D_i$ as

\begin{equation}\label{D}
\begin{aligned}
\hat D_i&=\int_0^\infty (\log s+\frac{1}{\mu_q} \log\frac{\mu_p}{\mu_q})s^{2\gamma_q+1}\sum_j a_{ij}e^{\bar U_j^k(s)}ds\\
&=\bar D_i+\frac{\bar m_{iv}^k}{\mu_q} \log\frac{\mu_p}{\mu_q}
\end{aligned}
\end{equation}

By the dependence of initial condition, we have
\begin{equation}\label{D-comp-2}
|D_i-\hat D_i|\leq C\sum_{i=2}^{n}|\alpha_i-\bar\alpha_i|\leq C\sum_{i=2}^n\bigg|\frac{\sigma_i^k}{\mu_p}-\frac{\Bar \sigma_i^k}{\mu_q}\bigg|
\end{equation}
Combining (\ref{3.11}), (\ref{3.13}, (\ref{D}) and (\ref{D-comp-2}) we have

\begin{equation}\label{3.14}
\bigg|\bigg(\frac{\Bar m_{iv}^k}{2\mu_q}-\frac{m_{iv}^k}{2\mu_p}\bigg)M_{k,p} \bigg| \leq C\sum_{j=2}^n\bigg|\frac{\sigma_j^k}{\mu_p}-\frac{\Bar \sigma_j^k}{\mu_q}\bigg|+E
\end{equation}

After multiplying $a^{ij}$ with summation on $i$ and taking summation on $j$, we have
\begin{equation*}
\sum_{j=2}^n\bigg|\frac{\Bar\sigma_j^k}{\mu_q}-\frac{ \sigma_j^k}{\mu_p}\bigg|\leq \frac{C}{M_{k,p}}\sum_{l=2}^n\bigg|\frac{\sigma_l^k}{\mu_p}-\frac{\Bar \sigma_l^k}{\mu_q}\bigg|+E/M_{k,p}.
\end{equation*}
Hence we obtain 
\begin{equation}\label{3.15}
\frac{\sigma_i^k}{\mu_p}-\frac{\Bar\sigma_i^k}{\mu_q}=O(\varepsilon_{k,p}^{m_{k,p}-2\mu_p})/M_{k,p},\quad i\in I.
\end{equation}
and have proved that
$$\frac{m_{iv}^k}{\mu_p}-\frac{\Bar m_{iv}^k}{\mu_q}=O(\varepsilon_{k,p}^{m_{k,p}-2\mu_p})/M_{k,p},$$
using the expansion of $\tilde u_i^k$, we get
$$|m_i^k-m_{iv}^k|\leq C\delta_0^{2+2\mu_p-m_k}\varepsilon_{k,p}^{m_{k,p}-2\mu_p},\quad |\Bar m_i^k-\Bar m_{iv}^k|\leq C\delta_0^{2+2\mu_p-m_k}\varepsilon_{k,p}^{m_{k,p}-2\mu_p}$$

Thus
$$\bigg|\frac{m_{i}^k}{\mu_p}-\frac{\Bar m_{i}^k}{\mu_q}\bigg|\leq C\delta_0^{2+2\mu_p-m_k}\varepsilon_{k,p}^{m_{k,p}-2\mu_p}$$

Proposition \ref{close-m} is established. $\Box$

\begin{rem}Proposition \ref{close-m} can be easily applied to multiple-blowup-point cases: If there are $L_2$ blowup points $p_1,...,p_{L_2}$, we have
$$|\frac{m_{it}^k}{\mu_t}-\frac{m_{is}^k}{\mu_s}|\le C\delta_0^{2+\mu_t-m_{t,k}}\epsilon_{k,t}^{m_{k,t}-2\mu_t}, \quad t,s=1,...,N, $$
if $\lim_{k\to \infty}\frac{m_{k,t}}{\mu_t}<4$ or $\lim_{k\to \infty}m_{k,t}/\mu_t=4$ but $\gamma_t,\gamma_s\neq 0$. 
\end{rem}

\medskip

As a consequence of Proposition \ref{close-m} we can obtain a rather accurate estimate of $\rho_i^k-\rho$ at this stage.  If we write $\rho_i^k$ as $\int_{M}\rho_i^k\mathfrak{h}_i^ke^{\hat u_i^k}dV_g$, we can further write it as
$$\rho_i^k=\sum_{t=1}^{L_2}\rho_{it}^k+\rho_{ib}^k$$
where 
$$\rho_{it}^k=\int_{B(p_t,\delta)}\rho_i^k\mathfrak{h}_i^ke^{\hat u_i^k} \quad \mbox{and}\quad 
\rho_{ib}^k=\int_{M\setminus \cup_{t}B(p_t,\delta)}\rho_i^k\mathfrak{h}_i^ke^{\hat u_i^k}.$$
A trivial estimate on $\rho_{ib}^k$ is $\rho_{ib}^k=E$ since outside bubbling disks, $\hat u_i^k=\bar u_i^k+O(1)$ and $e^{\bar u_i^k}=E$. When the equation is written in local coordinates around $p_t$, we use the notation 
$$\sigma_{it}^k:=\frac 1{2\pi}\int_{B(p_t,\delta)}|y|^{2\gamma_t}\tilde h_i^ke^{\tilde u_i^k}dV_g, \quad m_{it}^k=\sum_{j=1}^na_{ij}\sigma_{it}^k. $$
Then clearly $\sigma_{it}^k=\rho_{it}^k/2\pi$. 
The conclusion of Proposition \ref{close-m} gives
$$\rho_{it}^k/\mu_t-\rho_{is}^k/\mu_s=E, $$
Next we claim that 
\begin{equation}\label{close-m-3}
\frac{m_{it}^k}{\mu_t}-\mathfrak{m}_i=E
\end{equation}

This can be verified from Pohozaev identities. For each $t$, we have
$$\sum_i 4\mu_t\sigma_{it}^k-\sum_{i,j}a_{ij}\sigma_{it}^k\sigma_{jt}^k=O(\epsilon_{t,k}^{m_t^k-2\mu_t}).$$
Taking the sum of $t$ we have 
$$\sum_i 4\frac{\rho_i^k}{2\pi n_L}-\sum_{ij}a_{ij}\frac{\rho_i^k}{2\pi n_L}\frac{\rho_j^k}{2\pi n_L}=E. $$
If we write $\rho_i=\rho_{i}^k+s_i^k$, then the difference between the equations for $\rho^k$ and $\rho$ gives
$$2\sum_i(2-\sum_j a_{ij}\frac{\rho_j^k}{2\pi n_L})s_i^k=E. $$
since $\sum_j a_{ij}\frac{\rho_j^k}{2\pi n_L}>2$ for each $i$ and all $s_i^k$s have the same sign (see the assumption 
(\ref{rho-com}), we have deduced that $\rho_i-\rho_i^k=E$ for all $i$. 

By the definition of $\mathfrak{m}$ in (\ref{def-m}) and the fact that 
$\sum_t\mu_t=n_L$ (\ref{relation-1}) we see that $E$ can be written as $O(\epsilon_k^{\mathfrak{m}-2})$ if either $\mathfrak{m}<4$ or $\mathfrak{m}=4$ and all the blowup points are singular sources. If $\mathfrak{m}=4$ there exists a regular blowup point, $E=O(\epsilon_k^2\log \frac{1}{\epsilon_k})$. 

We also use $V_t^k=(V_{1t}^k,\cdots,V_{nt}^k)$ to denote the sequence of global solutions that approximate $\tilde u_i^k$ in $B(p_t^k,\delta_0)$. In the context of multiple bubbles, $V_{it}^k$ has two expressions: First from the evaluation of $\tilde u_i^k$ away from blowup points (obtained from the Green's representation of $\hat u_i^k$) we have 
\begin{equation}\label{3.17}
\begin{aligned}
V_{it}^k(x)&=-m_{it}^k\log|x|+\Bar{u}_i^k+2\pi \sum_{l=1}^{L_2} m_{il}^k G^*(p_t^k,p_l^k)+\log(\rho_i^k h_i^k(p_t^k))+E
\end{aligned}
\end{equation}
and, on the other hand, from the expansion of $U_i^k$ and scaling to $V_i^k$ we have
\begin{equation}\label{3.18}
\begin{aligned}
V_{it}^k(x)=-m_{itv}^k\log|x|-\frac{m_{itv}^k-2\mu_t}{2}M_{kt}+D_{it}-\alpha_{it}+E.
\end{aligned}
\end{equation}

By comparing the two expressions of $V_{it}^k$ we have
\begin{equation*}
\begin{aligned}
-m_{it}^k\log|x|+\Bar{u}_i^k+2\pi \sum_{l=1}^{L_2} m_{il}^k G^*(p_t^k,p_l^k)+\log(\rho_i^k h_i^k(p_t^k))\\
=-m_{itv}^k\log|x|-\frac{m_{itv}^k-2\mu_t}{2}M_{k,t}+D_{it}-\alpha_{it}+E.
\end{aligned}
\end{equation*}
for $x$ outside bubbling disks. 
Solving $\Bar u_i^k$ from the above, we have
\begin{equation}\label{3.19}
\begin{aligned}
\Bar u_i^k=&-\frac{m_{it}^k-2\mu_t}{2}M_{k,t}-2\pi \sum_{l=1}^{L_2} m_{il}^k G^*(p_t^k,p_l^k)-\log(\rho_i^k h_i^k(p_t^k))\\
&+D_{it}-\alpha_{it}+E.
\end{aligned}
\end{equation}

Thus, for $t \neq s$, we observe that 
$$|(m_{it}^k-2\mu_t)-(m_{is}^k-2\mu_s)|\le C\epsilon_k^{\mathfrak{m}-2}/M_k, $$ 
$$|D_{it}-D_{is}-\mathfrak{m}_i\log \frac{\mu_t}{\mu_s}|+|\alpha_{it}-\alpha_{is}|=E, $$
where we have used 
\begin{equation}\label{d-relation}
D_{it}=D_{is}+\mathfrak{m}_i\log \frac{\mu_t}{\mu_s}+E.
\end{equation}
which comes from (\ref{D}). If $\rho=Q$ and there is a regular blowup point, the error is changed to $O(\epsilon_k^2\log 1/\epsilon_k)$.
Thus 
\begin{align}\label{bar-u-i}
\bar u_i^k=-(\mathfrak{m}_i-2)\frac{\mu_t M_{k,t}}{2}-2\pi\mathfrak{m_i}
\sum_{l=1}^{L_2} \mu_l G^*(p_t,p_l)\\
-\log (\rho_i^kh_i^k(p_t^k))+D_i-\alpha_i+\mathfrak{m}_i\log \frac{\mu_t}{\mu_1}
+E. \nonumber
\end{align}
As a consequence 
\begin{equation}\label{exp-bar-u}
e^{\bar u_i^k}=\epsilon_{k,t}^{\mu_t(\mathfrak{m}_i-2)}e^{-2\pi \mathfrak{m}_i\sum_l\mu_l
G^*(p_t^k,p_l^k)}e^{D_{i}-\alpha_{i}}(\frac{\mu_t}{\mu_1})^{\mathfrak{m}_i}/(\rho_i^kh_i^k(p_t^k))+O(\epsilon_k^{\mathfrak{m}-2+\delta})
\end{equation}
where $\delta>0$ is independent of $k$, $D_i=D_{i1}$, $\alpha_i=\alpha_{i1}$ and (\ref{d-relation}) is used.

Comparing the two different expressions of $\bar u_i^k$ in terms of $t$ and $s$, we have 
\begin{align}\label{rel-2}
&(\mathfrak{m}_i-2)\frac{M_{t,k}\mu_t-M_{s,k}\mu_s}2
+2\pi \mathfrak{m}_i(\sum_{l=1}^{L_2}\mu_l(G^*(p_t^k,p_l^k)-G^*(p_s^k,p_l^k))\\
&+\log \frac{h_i^k(p_t^k)}{h_i^k(p_s^k)}+\mathfrak{m}_i\log \frac{\mu_s}{\mu_t}=E.\nonumber 
\end{align}

(\ref{rel-2}) can be used to simplify this notation:
\begin{align}\label{comp-e}
&\epsilon_{kt}^{\mu_t(\mathfrak{m}_i-2)}/\epsilon_{ks}^{\mu_s(\mathfrak{m}_i-2)}\\
=&e^{2\pi \mathfrak{m}_i(\sum_l \mu_l (G^*(p_t^k,p_l^k)-G^*(p_s^k,p_l^k))}
\frac{h_i^k(p_t^k)}{h_i^k(p_s^k)}(\frac{\mu_s}{\mu_t})^{\mathfrak{m}_i}+O(\epsilon_k^{\delta}). \nonumber
\end{align} 
for some $\delta>0$. 

Another consequence of (\ref{rel-2}) is Theorem \ref{relation-h}:

\noindent{\bf Proof of Theorem \ref{relation-h}:}  From (\ref{rel-2}) we have 
\begin{equation}\label{add-1}
    \frac{M_{t,k}\mu_t-M_{s,k}\mu_s}2+H_{it}^k-H_{is}^k=E 
\end{equation}
where 
\begin{equation}\label{Hik-1}
H_{it}=\frac{2\pi\mathfrak{m}_i}{\mathfrak{m}_i-2}(\sum_{l=1}^{L_2}\mu_l G^*(p_t^k,p_l^k))+
\frac{1}{\mathfrak{m}_i-2}\log \frac{h_i^k(p_t^k)}{\mu_t^{\mathfrak{m}_i}}.
\end{equation}
Since the first term of (\ref{rel-2}) and the left hand side of (\ref{add-1}) are independent of $i$, we have 
\begin{equation}\label{Hik-2}
H_{it}-H_{is}=H_{jt}-H_{js}+E, \quad \forall i,j=1,...,n, \quad t,s=1,...,L_2. 
\end{equation}
Theorem \ref{relation-h} is established. $\Box$

\medskip

\noindent{\bf Proof of Theorem \ref{localmass}:}
Theorem \ref{localmass} follows directly from Proposition \ref{close-m} and (\ref{close-m-3}). $\Box$
\section{The leading terms in approximations }

\noindent{\bf Proof of Theorem \ref{gammaNpneq}:} First we recall that $\int_M\mathfrak{h}_i^ke^{\hat u_i^k}=1$. Then we write $\rho_i^k$ as
$$\rho_i^k=\sum_{t=1}^{L_2}\int_{B(p_t^k,\delta_0)} \rho_i^k \mathfrak{h}_i^k e^{\hat u_i^{k}}dV_g+\int_{M\setminus\cup_t B(p_t^k,\delta_0)}\rho_i^k \mathfrak{h}_i^k e^{\hat u_i^{k}}dV_g=\sum_{t=1}^{L_2}\rho_{it}^k+\rho_{ib}^k, $$
The notation $I_1$ in the introduction is 
$$I_1=\{i\in I; \mathfrak{m}_i=\mathfrak{m} \}$$

If we use $V_i^k$ to be the leading term in the approximation of $\tilde u_i^k$ and $U_i^k$ be the scaled version of $V_i^k$, by (\ref{16}) we have

\begin{equation}\label{3.21}
\frac{1}{2\pi \mu_t}\rho_{it}^k=\frac{\sigma_{it}^k}{\mu_t}-\frac{e^{D_{it}-\alpha_{it}}}{\mu_t(m_t-2\mu_t)}\varepsilon_{k,t}^{\kappa_t}\delta_0^{-\kappa_t}+E_{\delta_0}
\end{equation}
where $\sigma_{it}^k$ is the total integration of the approximating global solutions around $p_t^k$, 
$\kappa_t=m_{t,k}-2\mu_t$. 
Now for $i \notin I_1$ we have
\begin{equation}\label{3.22}
\frac{1}{2\pi \mu_t}\rho_{it}^k=\frac{\sigma_{it}^k}{\mu_t}+E_{\delta_0},\quad i \notin I_1
\end{equation}
and
\begin{equation}\label{3.23}
|\rho_{ib}^k|=E_{\delta_0}, \quad i \notin I_1
\end{equation}
For $t\neq s$ we have 
$$\frac{\sigma_{it}^k}{\mu_t}=\frac{\sigma_{is}^k}{\mu_s}+O(\epsilon_k^{\mathfrak{m}-2})/M_{k}.$$
From Remark \ref{pohozaev-R} we have 
$$
 4\sum_{i=1}^n \frac{\rho_{it}^k}{2\pi\mu_t}-\sum_{i=1}^n\sum_{j=1}^n a_{ij} \frac{\rho_{it}^k}{2\pi\mu_t}\frac{\rho_{jt}^k}{2\pi\mu_t} 
=\frac{2}{\mu_t^2} \delta_0^{-\kappa_t}\varepsilon_{k,t}^{\kappa_t}\sum_{i\in I_1}e^{D_{it}-\alpha_{it}}+E_{\delta_0}.
$$
Here we recall that for $x$ away from bubbling disks, 
\begin{align*}
\hat u_i^k(x)=\bar u_i^k+\sum_{t=1}^{L_2}2\pi m_{it}^kG(x,p_t^k)+E\\
=\bar u_i^k+\sum_{t=1}^{L_2}2\pi \mathfrak{m}_i\sum_{t=1}^{L_2}\mu_tG(x,p_t^k)+E. 
\end{align*}
Now we can obtain a more precise expression of $\rho_{ib}^k$ from the estimate of $\bar u_i^k$ in (\ref{3.19}):
\begin{flalign}\label{3.24}
&\rho_{ib}^k=\int_{M\setminus \cup_t B(p_t,\delta_0)} \rho_i^k \mathfrak{h}_i^k e^{\hat u_i^k}dV_g&\\
&=\int_{M\setminus \cup_t B(p_t^k,\delta_0)} \rho_i^k \mathfrak{h}_i^k e^{\bar u_i^k}e^{2\pi \mathfrak{m}_i\sum_{l=1}^{L_2} \mu_l G(x,p_l^k)}+E_{\delta_0} \nonumber &\\
&=\sum_{t=1}^{L_2}\varepsilon_{k,t}^{\kappa_t}\int_{\hat \Omega_{t,\delta}} \frac{|x-p_t^k|^{2\gamma_t}h_i^k(x)}{h_i^k(p_t^k)}e^{D_{it}-\alpha_{it}}e^{2\pi \mathfrak{m}_i\sum_{l=1}^{L_2} \mu_l(G(x,p_l^k)-G^*(p_t^k,p_l^k))}dV_g+E_{\delta_0} \nonumber
\end{flalign}

Now the leading term of $\Lambda_{L}$ can be simplified:
\begin{align*}
    &4\sum_i \frac{\rho_i^k}{2\pi n_L}-\sum_{i,j}a_{ij}\frac{\rho_i^k}{2\pi n_L}\frac{\rho_j^k}{2\pi n_L}\\
    =&\frac{1}{4\pi^2 n_L^2}(8\pi n_L \sum_i \rho_i^k-\sum_{i,j}a_{ij}\rho_i^k\rho_j^k)\\
    =&\frac{1}{4\pi^2n_L^2}(16\pi^2n_L\sum_i \frac{\rho_i^k}{2\pi}-4\pi^2\sum_{ij}a_{ij}
    \frac{\rho_i^k}{2\pi}\frac{\rho_j^k}{2\pi})\\
    =&\frac{1}{4\pi^2 n_L^2}\bigg (16\pi^2 n_L\sum_i (\sum_t \frac{\sigma_{it}^k}{\mu_t}\mu_t+\frac{\rho_{ib}^k}{2\pi})-4\pi^2\sum_{ij}a_{ij}
    (\sum_t\frac{\sigma_{it}^k}{\mu_t}\mu_t+\frac{\rho_{ib}^k}{2\pi})(\sum_s \frac{\sigma_{js}^k}{\mu_s} \mu_s+\frac{\rho_{jb}^k}{2\pi})\bigg )
\end{align*}
where the insignificant error is ignored. Using $\frac{\sigma_{it}^k}{\mu_t}=\frac{\sigma_{is}^k}{\mu_s}+o(\epsilon_k^{\mathfrak{m}-2})$ we have
\begin{align*}
    &4\sum_i \frac{\rho_i^k}{2\pi n_L}-\sum_{i,j}a_{ij}\frac{\rho_i^k}{2\pi n_L}\frac{\rho_j^k}{2\pi n_L}\\
    =&\frac{4}{n_L}\sum_i\sum_t(\frac{\sigma_{it}^k}{\mu_t}\mu_t+\frac{\rho_{ib}^k}{2\pi})
    -\frac{1}{n_L^2}\sum_{ij}a_{ij}(\sum_t\frac{\sigma_{it}^k}{\mu_t}\mu_t)(\sum_s\frac{\sigma_{is}^k}{\mu_s}\mu_s)-\frac{2}{n_L^2}\sum_{ij}a_{ij}(\sum_t\frac{\sigma_{it}^k}{\mu_t}\mu_t)\frac{\rho_{jb}^k}{2\pi}\\
    =&\frac{4}{n_L}\sum_i(\sum_t \frac{\sigma_{it}^k}{\mu_t}\mu_t+\frac{\rho_{ib}^k}{2\pi})-\frac{1}{n_L^2}\sum_{ij}a_{ij}(\sum_t \mu_t
    \frac{\sigma_{it}^k}{\mu_t}\frac{\sigma_{jt}^k}{\mu_t})n_L-\frac{2}{n_L}\sum_i\mathfrak{m}_i\frac{\rho_{ib}^k}{2\pi}\\
    =&\frac{1}{n_L}\sum_t\frac{2}{\mu_t}\delta_0^{\mu_t(2-\mathfrak{m}_i)}\epsilon_{k,t}^{\mu_t(\mathfrak{m}_i-2)}\sum_{i\in I_1}e^{D_{it}-\alpha_{it}}-\frac{2}{n_L}\sum_i(\mathfrak{m}_i-2)\frac{\rho_{ib}^k}{2\pi}. 
    \end{align*}

\begin{align*}
    \rho_{ib}^k=&\sum_{t=1}^{L_2}\int_{\hat \Omega_{t,\delta_0}} \rho_i^k\mathfrak{h}_i^ke^{\hat u_i^k}dV_g
    =\sum_{t=1}^{L_2}\int_{\hat \Omega_{t,\delta}} \rho_i^k\mathfrak{h}_i^ke^{\bar u_i^k+2\pi\mathfrak{m}_i\sum_l \mu_l G(x,p_l^k)}\\
    =&\sum_{t=1}^{L_2}\int_{\hat \Omega_{t,\delta}} \rho_i^k\mathfrak{h}_i^ke^{\bar u_i^k+2\pi\mathfrak{m}_i(-\frac{\mu_t}{2\pi}\log |x-p_t^k|+\sum_l \mu_l G^*(x,p_l^k))}\\
    =&\sum_{t=1}^{L_2}\int_{\hat \Omega_{t,\delta}} \rho_i^k h_i^k(x)e^{\bar u_i^k}|x-p_t^k|^{2\gamma_t-\mu_t\mathfrak{m}_i}e^{2\pi\mathfrak{m}_i\sum_l\mu_l G^*(x,p_l^k)}
    \end{align*}
Using (\ref{3.19}) we write $\rho_{ib}^k$ as
$$\rho_{ib}^k=\sum_t\int_{\hat\Omega_{t,\delta_0}}\frac{h_i^k}{h_i^k(p_t^k)}\epsilon_{t,k}^{\kappa_t}e^{D_{it}-\alpha_{it}}|x-p_t^k|^{2\gamma_t-\mu_t\mathfrak{m}_i}e^{2\pi\mathfrak{m}_i\sum_l(\mu_l(G^*(x,p_l^k)
-G^*(p_t^k,p_l^k))}+E.$$

Thus
\begin{align*}
    &4\sum_i\frac{\rho_i^k}{2\pi n_L}-\sum_{ij}a_{ij}\frac{\rho_i^k}{2\pi n_L}\frac{\rho_j^k}{2\pi n_L}
    =\frac{2}{n_L}\sum_{i}\sum_{t=1}^{L_2}\epsilon_{k,t}^{\kappa_t}e^{D_{it}-\alpha_{it}}
    (\frac{\delta_0^{\mu_t(2-\mathfrak{m}_i)}}{\mu_t} \\
   & -\int_{\hat \Omega_{t,\delta_0}}(\frac{\mathfrak{m}_i-2}{2\pi})\frac{h_i^k(x)}{h_i^k(p_t^k)}|x-p_t^k|^{2\gamma_t-\mu_t\mathfrak{m}_i}e^{2\pi \mathfrak{m}_i\sum_{l=1}^{L_2}\mu_l(G^*(x,p_l^k)-G^*(p_t^k,p_l^k)}).
\end{align*}
From here we see that when $\delta_0\to 0$, the leading term $O(\delta_0^{\mu_t(2-\mathfrak{m}_i)})$ is cancelled out. Using the relation between $\epsilon_{t,k}$ and $\epsilon_k$ in (\ref{comp-e}), $D_{it}$ and $D_i$ in (\ref{d-relation}) $\alpha_{it}=\alpha_t+E_{\delta}$ we have
\begin{align*}
 &4\sum_i\frac{\rho_i^k}{2\pi n_L}-\sum_{ij}a_{ij}\frac{\rho_i^k}{2\pi n_L}\frac{\rho_j^k}{2\pi n_L}\\
 =&\frac{2}{n_L}\epsilon_k^{\mathfrak{m}-2}\sum_{i\in I_1}\sum_t B_{it}
 (\frac{\delta_0^{\mu_t(2-\mathfrak{m})}}{\mu_t}\\
 &-\int_{\hat \Omega_{t,\delta_0}}(\frac{\mathfrak{m}-2}{2\pi})\frac{h_i^k(x)}{h_i^k(p_t^k)}|x-p_t^k|^{2\gamma_t-\mu_t\mathfrak{m}}e^{2\pi \mathfrak{m}\sum_{l=1}^{L_2}\mu_l(G^*(x,p_l^k)-G^*(p_t^k,p_l^k)})\\
 =&\frac{2}{n_L}\epsilon_k^{\mathfrak{m}-2}\sum_{i\in I_1}\sum_t B_{it}
 (\frac{\delta_0^{2\mu_t-m_t}}{\mu_t}\\
 &-\int_{\hat \Omega_{t,\delta_0}}(\frac{\mathfrak{m}-2}{2\pi})\frac{h_i^k(x)}{h_i^k(p_t^k)}|x-p_t^k|^{(2-\mathfrak{m})\mu_t-2}e^{2\pi \mathfrak{m}\sum_{l=1}^{L_2}\mu_l(G^*(x,p_l^k)-G^*(p_t^k,p_l^k)})
\end{align*}
where 
$$B_{it}=e^{2\pi\mathfrak{m}(\sum_{l=1}^{L_2}\mu_l(G^*(p_t^k,p_l^k)-G^*(p_1^k,p_l^k))}\frac{h_i^k(p_t^k)}{h_i^k(p_1^k)}e^{D_i-\alpha_i}
$$ 
if $m_{it}=m_t$, otherwise $B_{it}$ is any constant.  Here $D_i=D_{i1}$ and $\alpha_i=\alpha_{i1}$. Theorem \ref{gammaNpneq} is established. $\Box$

\bigskip

\noindent{\bf Proof of Theorem \ref{main-3}:}  Around each regular blowup point $p_t^k$, we have
$$\sum_{i=1}^n(\nabla (\log h_i^k+\phi_i^k)(p_t^k))\sigma_{it}^k=O(\epsilon_k^{\mathfrak{m}-2}), \quad \mbox{if}\quad m_k<4\mu_p,$$
and 
$$\sum_{i=1}^n(\nabla (\log h_i^k+\phi_i^k)(p_t^k))\sigma_{it}^k=O(\epsilon_k^2\log \frac{1}{\epsilon_k}), \quad \mbox{if}\quad m_k=4\mu_p,$$
In the first case, 
$$\nabla \phi_i^k(p_t^k)=\sum_{l=1}^{L_2}2\pi m_l^k\nabla_1 G^*(p_t^k,p_l^k)+O(\epsilon_k^{\mathfrak{m}-2}) $$
Using $\frac{m_l^k}{\mu_l}=\mathfrak{m}+O(\epsilon_k^{\mathfrak{m}-2})$ we obtain 
$$\nabla \phi_i^k(p_t^k)=2\pi \mathfrak{m}\sum_{l=1}^{L_2}\mu_l \nabla_1 G^*(p_t^k,p_l^k))+O(\epsilon_k^{\mathfrak{m}-2}).$$
Since $\sigma_{i,t}^k=\frac{\rho_i}{2\pi n_L}+O(\epsilon_k^{\mathfrak{m}-2})$,
(\ref{11july13e7}) follows immediately. The derivation of (\ref{11july13e8}) can be derived in a similar fashion. Theorem \ref{main-3} is established. $\Box$

\medskip

\noindent{\bf Proof of Theorem \ref{m-eq-4}:}

$$\rho_i^k=\sum_{t=1}^{L_2}\int_{B(p_t^k,\delta_0)}\rho_i^k h_i^k e^{u_i^k} dV_g +\int_{M\setminus \cup_{t=1}^{L_2} B(p_t^k,\delta_0)} \rho_i^k h_i^k e^{u_i^k} dV_g. $$
We continue to use the notation $\rho_{it}^k$ and $\rho_{ib}^k$.
In this case the $\rho_{ib}^k=O(\epsilon_k^2)$, which is an error. The leading term comes from interior integration. 

Now we use the expansion of bubbles to compute each $\rho_{it}^k$. By the expansion of $\tilde u_i^k$ around $p_t^k$, which is a regular blowup point, we have

\begin{equation}\label{11july5e1}
\begin{split}
\rho_{it}^k=& \int_{B(p_t^k,\delta_0)}\rho_i^k h_i^k e^{u_i^k} dV_g
=\int_{B(0,\delta_0)}\tilde h_i^k e^{\phi_i^k} e^{\tilde u_i^k-\phi_i^k}d\eta \\
=&\int_{B(0,\delta_0 \epsilon_{k,t}^{-1})} \rho_i^k h_i^k(p_t^k) e^{U_i^k(\eta)}d\eta+O(\epsilon_{k,t}^2)\\
&+\frac 14\int_{B(0,\delta_0 \epsilon_{k,t}^{-1})}
\epsilon_{k,t}^2\big( \Delta (\log \tilde h_i^k)(0)+\frac 14|\nabla (\log \tilde h_i^k+\phi_i^k)(0)|^2)
|y|^2e^{U_i^k}d\eta. 
\end{split}
\end{equation}
Note that in the expansion of 
$\tilde u_i^k-\phi_i^k+(\log \tilde h_i^k+\phi_i^k)$, the first term is a global solution $V_k$, the second term is a projection onto $e^{i\theta}$ that leads to integration zero. The third term has the leading term $\Delta (\log \tilde h_i^k)(0)$
and the square of $\nabla (\log \tilde h_i^k+\phi_i^k)$. 

On the other hand if $p_t$ is a singular point, 
$$\rho_{it}^k=\int_{B(0,\delta_0 \epsilon_{k,t}^{-1})} \rho_i^k|y|^{2\gamma_t}
e^{U_i^k(y)}dy+O(\epsilon_k^2). $$

From (\ref{close-m-3}) we see that $\epsilon_{k,t}$ can be replaced by $\epsilon_k$. Hence the first integral on the right hand side of the above is $O(\epsilon_k^2)$ different from the global solution in the approximation of $\tilde u_i^k$ around $p_t^k$.
So we use $\sigma_{ivt}^k$ to denote it. For $t\neq s$, from (\ref{3.15}) we see that
$$\sigma_{ivt}^k-\sigma_{ivs}^k=O(\epsilon_k^2/\log \frac{1}{\epsilon_k}). $$
 To evaluate the last term, we first use the
definition of the $\tilde h_i^k$ to have

$$\Delta (\log \tilde h_i^k)(0)=\Delta (\log h_i^k)(0)-2K(p_t)+8\pi n_L+
O(\epsilon_k^2). $$
$$\nabla \log \tilde h_i^k(0)=\nabla \log h_i^k(p_k), $$
$$\nabla \phi_i^k(0)=8\pi \sum_{l=1}^{L_2}\mu_l \nabla_1 G^*(p_t^k,p_l^k)
+O(\epsilon_k^2\log \frac{1}{\epsilon_k}). $$
Then we define $b_{it}^k$ as
\begin{align}\label{bit}
b_{it}^k&=e^{D_i-\alpha_i}\bigg (\frac 14 \Delta (\log h_i^k)(p_t^k)-\frac{K(p_t^k)}2+2\pi n_L  \\
&+\frac 14|\nabla (\log h_i^k)(p_t^k)+8\pi\sum_{l=1}^{L_2}\nabla_1G^*(p_t^k,p_l^k)|^2
\bigg ) \nonumber
\end{align}
With this $b_{it}^k$ we have
$$\frac{\rho_{it}^k}{2\pi}=\sigma_{iv}^k+b_{it}^k\epsilon_{k}^2\log \frac{1}{\epsilon_k}+O(\epsilon_k^2). $$
Consequently,
\begin{equation}\label{11jun8e1}
\begin{split}
&4\sum_{i=1}^n\frac{\rho_i^k}{2\pi n_L}- \sum_{i=1}^n\sum_{j=1}^na_{ij}\frac{\rho_i^k}{2\pi n_L}\frac{\rho_j^k}{2\pi n_L}\\
=&4\sum_{i=1}^n\sum_{t=1}^{L_2}\frac{\rho_{it}^k}{2n_L \pi}- \sum_{i=1}^n\sum_{j=1}^na_{ij}(\sum_{t=1}^{L_2} \frac{\rho_{it}}{2\pi n_L})(\sum_{s=1}^{L_2} \frac{\rho_{js}}{2\pi n_L})  \\
=&4\sum_{i=1}^n\sum_{t=1}^{L_2}(\frac{\sigma_{ivt}}{n_L}+\tilde \epsilon_k\frac{b_{it}^k}{n_L})-\sum_{i=1}^n\sum_{j=1}^n\sum_{t=1}^{L_2}\sum_{s=1}^{L_2}a_{ij}(\frac{\sigma_{ivt}}{n_L}+\tilde\epsilon_k\frac{b^k_{it}}{n_L})(\frac{\sigma_{jv}}{n_L}+\tilde\epsilon_k\frac{b^k_{js}}{n_L})  \\
=&4\sum_{i=1}^n\frac{\sigma_{iv}}{n_L}- \sum_{i=1}^n\sum_{j=1}^na_{ij}\frac{\sigma_{iv}\sigma_{jv}}{n_L}+4\sum_{i=1}^n\sum_{t=1}^{L_2}\tilde\epsilon_k\frac{b^k_{it}}{n_L}
-2\sum_{i=1}^n\sum_{j=1}^n\sum_{t=1}^{L_2}a_{ij}\sigma_{iv}\tilde\epsilon_k\frac{b^k_{jt}}{n_L}\\
&+O(\epsilon_k^2)
 \\
=&-4\epsilon_k^2\log \epsilon_k^{-1}\sum_{i=1}^n\sum_{t=1}^{L_2} b_{it}^k+O(\epsilon_k^2).
\end{split}
\end{equation}
where $\tilde \epsilon_k$ stands for $\epsilon_k^2\log \frac{1}{\epsilon_k}$.
Theorem \ref{m-eq-4} is established. $\Box$

\section{Local approximation of bubbling solutions}  

In this section we provide asymptotic analysis for bubbling solutions with no oscillation on the boundary of a unit disk. The estimates of this section have been used repeatedly in the proof of the main theorems. Let $\mathfrak{u}^k=(\mathfrak{u}_1^k,...,\mathfrak{u}^k_n)$ be a sequence of blowup solutions of
\begin{equation}\label{uik}
-\Delta \mathfrak{u}_i^k=\sum_{j=1}^na_{ij}|x|^{2\gamma}\mathfrak{h}_j^ke^{\mathfrak{u}_j^k},\quad i=1,..,n,\quad x\in B_1,
\end{equation}
where $-1<\gamma<0$, 
\begin{equation}\label{hikn}
\frac 1C\le \mathfrak{h}_i^k(x)\le C,\quad \|\mathfrak{h}_i^k\|_{C^3(B_1)}\le C, \quad x\in B_1, \quad i=1,..,n.
\end{equation}
and the origin is the only blowup point in $B_1$:
$$\max_K\mathfrak{u}_i^k\le C(K), \quad \forall K\subset\subset \bar B_1\setminus \{0\}, \,\, \mbox{ and } \max_{B_1}\mathfrak{u}_i^k\to \infty,$$
with no oscillation on $\partial B_1$:
\begin{equation}\label{88e1}
\mathfrak{u}_i^k(x)-\mathfrak{u}_i^k(y)=0,\forall x,y\in \partial B_1,\quad i=1,...,n. 
\end{equation}
and uniformly bounded energy:
\begin{equation}\label{89e1}
\int_{B_1}|x|^{2\gamma}\mathfrak{h}_i^k e^{\mathfrak{u}_i^k}\le C,\quad C \mbox{ is independent of } k.
\end{equation}
Finally we assume that $\mathfrak{u}^k$ is a fully blown-up sequence, which means when re-scaled according its maximum, $\{\mathfrak{u}^k\}$ converges to a system of $n$ equations: Let
$$\bar M_k=\max_{B_1}\max_i\mathfrak{u}_i^k/\mu
\quad \mbox{and} \quad \bar\epsilon_k=e^{-\frac 1{2}\bar M_k},\quad \mu=1+\gamma$$

\subsection{First order estimates}\label{blowup}
Let $v^k=(v_1^k,...,v_n^k)$ be a scaling of $u^k$ according to its maximum:
\begin{equation}\label{vikdef}
v_i^k(y)=\mathfrak{u}_i^k(\bar \epsilon_ky)-2\mu\log \bar \epsilon_k,\quad y\in \Omega_k:=B(0,\bar \epsilon_k^{-1}).
\end{equation}
Then the equation for $v^k=(v_1^k,...,v_n^k)$ is
\begin{equation}\label{e-vk}
\Delta v_i^k+\sum_j a_{ij}|y|^{2\gamma}\mathfrak{h}^k_j(\bar \epsilon_k y)e^{v_j^k}=0,\quad
\mbox{in}\quad \Omega_k, \quad i=1,..,n,
\end{equation}
$v_i^k=\mbox{constant}$ on $\partial \Omega_k$ and
$v^k$ converges in $C^{\alpha}_{loc}(\mathbb R^2)$ for some $\alpha\in (0,1)$
to $U=(U_1,..,U_n)$, which satisfies
\begin{equation}\label{817e1}
\left\{\begin{array}{ll}-\Delta U_i=\sum_{j}a_{ij}|y|^{2\gamma}e^{U_j},\quad \mathbb R^2,\quad i=1,..,n\\
\\
\int_{\mathbb R^2}|y|^{2\gamma}e^{U_i}<\infty,\quad i=1,..,n,\quad \max_iU_i(0)=0.
\end{array}
\right.
\end{equation}
where for simplicity we assumed that $\lim_{k\to \infty}\mathfrak{h}_j^k(0)=1$. By a classification result of Lin-Zhang \cite{lin-zhang-dcds} that for $\gamma\in (-1,0)$
all the global solutions with finite energy are radial functions.
Set $V^k=(V_1^k,..,V_n^k)$ be the radial solutions of
\begin{equation}\label{88e2}
\left\{\begin{array}{ll}-\Delta V_i^k=\sum_{j=1}^na_{ij}|x|^{2\gamma}e^{V_j^k}\quad \mathbb R^2,\quad i\in I\\  \\
V_i^k(0)=\mathfrak{u}_i^k(0),\quad i\in I
\end{array}
\right.
\end{equation}
Here for simplicity we assume $\mathfrak{h}_i^k(0)=1$. Without this assumption the proof is still almost the same. 
 It is easy to see that any radial solution $V^k$ of (\ref{88e2}) exists for all $r>0$.
In \cite{linzhang1} the authors prove that
 \begin{equation}\label{11mar9e1}
 |\mathfrak{u}_i^k(x)-V_i^k(x)|\le C,\quad \mbox{for } |x|\le 1.
 \end{equation}
 From (\ref{11mar9e1}) we have the following spherical Harnack inequality:
 \begin{equation}\label{11mar9e3}
 |\mathfrak{u}_i^k(\bar \epsilon_ky)-\mathfrak{u}_i^k(\bar \epsilon_k y')|\le C
 \end{equation}
 for all $|y|=|y'|=r\le \bar \epsilon_k^{-1}$ and $C$ is a constant independent of $k,r$. (\ref{11mar9e3}) will play an essential role in the first order estimate.  Here we fix some notations:

$$\bar \sigma_i^k=\frac 1{2\pi}\int_{B_1}|x|^{2\gamma}\mathfrak{h}_i^ke^{\mathfrak{u}_i^k}, 
\quad \bar m_i^k=\sum_j a_{ij}\bar \sigma_j^k, $$
let $\bar \sigma_i$ and $\bar m_i$ be their limits and $\bar m=\min\{\bar m_1,...,\bar m_n\}$. Correspondingly we use $\sigma_{iv}^k$ and $m_{iv}^k$ to denote the energy for $V^k$:
$$\sigma_{iv}^k=\frac 1{2\pi}\int_{\mathbb R^2}|x|^{2\gamma}e^{V_i^k},\quad m_{iv}^k=\sum_j a_{ij}\sigma_{jv}^k. $$
Later we shall show that $\bar \sigma_i^k-\sigma_{iv}^k$ is very small. Right now we use the fact that
$\bar m\in (2\mu, 4\mu]$. Here is the reason: $\sigma^v=(\sigma_{1v},...,\sigma_{nv})$ satisfies
$$\sum_{ij}a_{ij}\sigma_{iv}\sigma_{jv}=4\mu\sum_i \sigma_{iv}..$$
This can be written as
$$\sum_i (m_{iv}-4\mu)\sigma_{iv}=0. $$
Thus $\min_i m_{iv}\le 4\mu$. Each $m_i>2\mu$ is by the integrability. 
Let $U^k=(U_1^k,...,U_n^k)$ be defined by 
\begin{equation}\label{def-Uk}
U_i^k(y)=V_i^k(\bar \epsilon_k y)+2\mu\log \bar \epsilon_k,\quad i\in I
\end{equation}
$U^k=(U_1^k,...,U_n^k)$ is the first term in the approximation to $v_i^k(y)$. Note that even though $U^k\to U$ we cannot use $U$ as the first term in the approximation because it does not have the same initial condition of $v^k$. 
 \begin{thm}\label{thm2}
Let $\delta\in (0,1)$ be defined by
$$\delta=\left\{\begin{array}{ll} 1-\epsilon, \quad \mbox{if } m-2\mu>1\\
1+2\mu-m+\epsilon, \quad \mbox{if}\quad m-2\mu\le 1. 
\end{array} 
\right.$$
Then there exists $C$ independent of $k$ such that 
\begin{align}\label{11mar9e2}
&|D^{\alpha}(v_i^k(y)-U_i^k(y)|\\
\le &
C\epsilon_k(1+|y|)^{\delta-|\alpha |}
\quad |y|<\bar \epsilon_k^{-1},\quad |\alpha|=0,1. \nonumber
\end{align}
\end{thm}

We use the following notations:
$$\bar \sigma_{iv}^k=\frac 1{2\pi}\int_{\mathbb R^2} |x|^{2\gamma} e^{V_i^k},\,\, \bar m_{iv}^k=\sum_{j=1}^n a_{ij}\bar \sigma_{jv}^k,\,\, \bar m_v^k=\min \{\bar m_{1v}^k,..,\bar m_{nv}^k\}. $$
From Theorem \ref{thm2} it is easy to see that $\lim_{k\to\infty}\bar m_i^k=\bar m_i$. Thus $\bar m_i^k\ge 2\mu +\delta_0$ for some $\delta_0>0$ independent of $k$. Let 
$$w_i^k(y)=v_i^k(y)-U_i^k(y),\quad i\in I. $$
 We first write the equation of $w_i^k$ based on (\ref{e-vk}) and (\ref{88e2}):
\begin{equation}\label{wik01}
\left\{\begin{array}{ll}
\Delta w_i^k(y)+\sum_ja_{ij}|y|^{2\gamma}\mathfrak{h}_j^k(\bar\epsilon_ky)e^{\xi_j^k}w_j^k=-\sum_ja_{ij}|y|^{2\gamma}(\mathfrak{h}_j^k(\epsilon_ky)-1)e^{U_j^k}, \\ \\
w_i^k(0)=0,\quad i\in I,
\end{array}
\right.
\end{equation}
where
$\xi_i^k$ is defined by
\begin{equation}\label{1211e1}
e^{\xi_i^k}=\int_0^1e^{tv_i^k+(1-t)U_i^k}dt.
\end{equation}

Since both $v^k$ and $U^k$ converge to $v$, $w_i^k=o(1)$ over any compact subset of $\mathbb R^2$.
The first estimate of $w_i^k$ is the following
\begin{lem}\label{1020lem1}
\begin{equation}\label{820e1}
w_i^k(y)=o(1)\log (1+|y|)+O(1),\quad \mbox{ for } y\in \Omega_k.
\end{equation}
\end{lem}

\medskip

\noindent{\bf Proof}: By (\ref{11mar9e3})
$$|v_i^k(y)-\bar v_i^k(|y|)|\le C, \quad \forall y\in \Omega_k $$
where $\bar v_i^k(r)$ is the average of $v_i^k$ on $\partial B_r$:
$$\bar v_i^k(r)=\frac{1}{2\pi r}\int_{\partial B_r}v_i^k. $$
 Thus we have
$|y|^{2\gamma} e^{v_i^k(y)}=O(r^{-2-\delta_0})$
and $|y|^{2\gamma}e^{U_i^k(y)}=O(r^{-2-\delta_0})$ where $r=|y|$ and $\delta_0>0$. Then
\begin{equation}\label{to-be-used}
r(\bar w_i^k)'(r)=\frac 1{2\pi }\bigg (\int_{B_r}\sum_ja_{ij}|y|^{2\gamma}\mathfrak{h}_j^k(\bar \epsilon_k y)e^{v_j^k}
-\int_{B_r}\sum_ja_{ij}|y|^{2\gamma}e^{U_j^k}\bigg )
\end{equation}

It is easy to use the decay rate of $e^{v_i^k}$, $e^{U_i^k}$ and the closeness between $v_i^k$ and $U_i^k$
to obtain
$$ r(\bar w_i^k)'(r)=o(1),\quad r\ge 1. $$
Hence $\bar w_i^k(r)=o(1)\log r$ and (\ref{820e1}) follows from this easily. Lemma \ref{1020lem1}
is established. $\Box$

\medskip

The following estimate is immediately implied by Lemma \ref{1020lem1}:
$$|y|^{2\gamma} e^{\xi_i^k(y)}\le C(1+|y|)^{-m+2\gamma+o(1)}\quad \mbox{ for } y\in \Omega_k.$$

Before we derive further estimate for $w_i^k$ we cite a useful estimate for the Green's function on $\Omega_k$ with respect to the
Dirichlet boundary condition (see \cite{lin-zhang-jfa}):

\begin{lem}(Lin-Zhang)\label{greenk} Let $G(y,\eta)$ be the Green's function with respect to Dirichlet boundary condition on $\Omega_k$. For $y\in \Omega_k$,
let \begin{eqnarray*}
\Sigma_1&=&\{\eta \in \Omega_k;\quad |\eta |<|y|/2 \quad \}\\
\Sigma_2&=&\{\eta \in \Omega_k;\quad |y-\eta |<|y|/2 \quad \}\\
\Sigma_3&=&\Omega_k\setminus (\Sigma_1\cup \Sigma_2).
\end{eqnarray*}
Then in addition for $|y|>2$,
\begin{equation}\label{1020e5}
|G(y,\eta)-G(0,\eta)|\le \left\{\begin{array}{ll}
C(\log |y|+|\log |\eta ||),\quad \eta\in \Sigma_1,\\
C(\log |y|+|\log |y-\eta ||),\quad \eta\in \Sigma_2,\\
C|y|/|\eta |,\quad \eta \in \Sigma_3.
\end{array}
\right.
\end{equation}
\end{lem}
Next we prove that there is no kernel in the linearized operator with controlled growth.

\begin{prop}\label{uniqlin} Let $U=(U_1,...U_n)$ be a solution of 
$$U_i''+\frac 1rU_i'(r)+\sum_j a_{ij}r^{2\gamma}e^{U_j}=0,\quad 0<r<\infty, \quad \int_0^{\infty}r^{2\gamma+1}e^{U_i}<\infty, \quad i=1,..,n,$$
where $-1<\gamma<0$. It is established in \cite{lin-zhang-dcds} that $U_i$ is radial there is a classification of all solutions $U$. 
Set $\delta_*=1-\epsilon$ for $\epsilon>0$ small.
Let $\phi=(\phi_1,...,\phi_n)$ be a solution of 
$$\Delta \phi_i+\sum_j |x|^{2\gamma} a_{ij}e^{U_j}\phi_j=0, \quad \mbox{in}\quad \mathbb R^2,\quad -1<\gamma<0, \quad i\in I$$
and 
$$\phi_i(x)=O((1+|x|)^{\delta_*})\quad i\in I.$$
If we further have $\phi_i(0)=0$ for $i\in I$, then $\phi_i\equiv 0$ for all $i$. 
\end{prop}

\noindent{\bf Proof of Proposition \ref{uniqlin}:} We are going to consider projections of $\phi$ on $e^{iN\theta}$ and prove that $\phi$ is actually bounded. Once this is proved, the conclusion is already derived in \cite{lin-zhang-dcds}. First we prove that the projection on $1$ is zero:
Let $\psi^0=(\psi_1^0,...,\psi_n^0)$ be the projection on $1$, then $\psi^0$ is the solution of 
$$\frac{d^2}{dr^2}\psi_i^0+\frac 1r\frac{d}{dr}\psi_i^0+\sum_j a_{ij}r^{2\gamma}e^{U_j}\psi_j^0=0, \quad 0<r<\infty, $$
with the initial condition $\psi_i^0(0)=0$. It is clear that as long as we prove $\psi_i^0\equiv 0$ in $[0,\tau]$ for a small
$\tau>0$, we are all done. Let 
$$M(r)=\max_{i\in I}\max |\psi_i^0(s)|,\quad 0\le s\le r. $$
Integration on the equation for $\psi_i^0$ with the zero initial condition gives 
$$\psi_i^0(r)=-\frac 1r\int_0^r\frac 1s\int_0^s\sum_j a_{ij}t^{1+2\gamma}e^{U_j(t)}\psi_j^0(t)dtds. $$
Using $M(r)$ for $\psi_j^0$ above, we have 
$$|\psi_i^0(r)|\le CM(r)r^{2\mu},\quad i\in I. $$
Obviously this leads to $M(r)\le Cr^{2\mu}M(r)$, which implies 
$M(r)=0$ for small $r$. Consequently $\psi_i^0\equiv 0$ for all $i\in I$. 

Next we consider the projection on $e^{il\theta}$ ( $l\ge 1$). 
 Let $\psi^l=(\psi^l_1,...,\psi^l_n)$ satisfy
$$\left\{\begin{array}{ll}
\frac{d^2}{dr^2}\psi^l_i(r)+\frac 1r\frac{d}{dr}\psi^l_i(r)-\frac {l^2}{r^2}\psi^l_i(r)=f_i, \quad 1\le i\le n, \quad \psi_i^l(0)=0, \\
\\
f_i(r)=-\sum_ja_{ij}r^{2\gamma}e^{U_j}\psi_j(r). 
\end{array} 
\right.
$$
solution of ode gives
$$\psi^l_i(r)=c_1r^l+c_2r^{-l}+r^l\int_{\infty}^r\frac{f(s)s^{1-l}}{2l}ds-r^{-l}\int_0^r\frac{s^{l+1}f(s)}{2l}ds$$
$\psi^l_i(0)=0$ gives $c_2=0$. $c_1=0$ because $\psi^l_i(r)$ has a sub-linear growth. Using $f(r)= O(r^{2\gamma-m+\delta_0})$ (where $\delta_0=1-\epsilon$), we obtain from standard evaluation that 
$$|\psi^l_i(r)|\le Cl^{-2}(1+r)^{2\mu-m+\delta_0}. $$
Taking the sum of all these projections we have obtained that 
$$\phi_i(x)=O((1+|x|)^{1-\epsilon-(m-2\mu)}+O(1). $$
The conclusion holds if $m-2\mu\ge 1$. If this is not the case, the same argument leads to 
$$\phi_i(x)=O((1+|x|)^{1-\epsilon-2(m-2\mu)}+O(1). $$
Obviously after finite steps we have $\phi_i(x)=O(1)$. By the uniqueness result of Lin-Zhang \cite{lin-zhang-dcds}, $\phi_i\equiv 0$. $\Box$

Proposition \ref{uniqlin} is established. $\Box$.

\bigskip

\begin{lem}\label{three-p-1}
Let 
$$\delta=\left\{\begin{array}{ll}
1-\epsilon,\quad \mbox{if}\quad m-2\mu\ge 1,\\
1+2\mu-m+\epsilon, \quad \mbox{if} \quad m-2\mu<1,
\end{array} 
\right. $$
there exists $C$ independent of $k$ such that 
$$|w_i^k|\leq C\bar \varepsilon_k (1+|y|)^\delta$$
\end{lem}

\noindent{\bf Proof of Lemma \ref{three-p-1}:} First we note that the $\delta$ defined above is less than $1$. 

Prove by contradiction, we assume

$$\Lambda_k:=\max_{y\in \Omega_k}\frac{\max_{i\in I}|w_i^k(y)|}{\bar \varepsilon_k(1+|y|)^\delta}\rightarrow\infty$$

Suppose $\Lambda_k$ is attained at $y_k\in \Bar{\Omega}_k$ for some $i_0\in I$.

$$\Bar{w}_i^k(y)=\frac{w_i^k(y)}{\Lambda_k\bar \varepsilon_k(1+|y_k|)^\delta}.$$

It follows from the definition of $\Lambda_k$ that 
\begin{equation}\label{lem3.18}
\Bar{w}_i^k(y)=\frac{|w_i^k(y)|}{\Lambda_k\bar \varepsilon_k(1+|y|)^\delta}\frac{(1+|y|)^\delta}{(1+|y_k|)^\delta}\leq \frac{(1+|y|)^\delta}{(1+|y_k|)^\delta}
\end{equation}

Then the equation for $\Bar{w}_i^k$ is

\begin{equation}\label{lem3.19}
  -\Delta \Bar{w}_i^k(y)=\sum_{j=1}^n a_{ij}|y|^{2\gamma}\mathfrak{h}_j^k(\bar \varepsilon_k y)e^{\xi_j^k}\Bar{w}_j^k(y)+\frac{C(1+|y|)^{1-m+2\gamma}}{\Lambda_k(1+|y_k|)^\delta}, \quad \Omega_k
\end{equation}

Here $\xi_i^k$ is given by (\ref{1211e1}). $\xi_i^k$ converges to $U_i$ in (\ref{v-global}) uniformly in any compact subset of $\mathbb R^2$.

First we claim that $|y_k|\to \infty$. Otherwise if there is a sub-sequence of $w_i^k$ (still denoted as $w_i^k$ that converges to $w=(w_1,...,w_n)$ of 
$$\Delta w_i+\sum_j a_{ij}|y|^{2\gamma}e^{U_j}w_j=0,\quad \mbox{in}\quad \mathbb R^2 $$
with $w_i(0)=0$ and $w_i(y)=O((1+|y|)^{\delta})$ for some $\delta<1$. Then $w_i\equiv 0$, which is a contradiction to $\max_i|w_i(y^*)|=1$ where $\lim_{k\to \infty} y_k=y^*$.

After ruling out $|y_k|\le C$, we now rule out $|y_k|\to \infty$ using the Green's representation formula of $\bar w_i^k$:

By Green's representation formula for $\Bar{w}_i^k$,

$$\Bar{w}_i^k(y)=\int_{\Omega_k}G(y,\eta)(-\Delta\Bar{w}_i^k(\eta))d\eta+\Bar{w}_i^k|_{\partial\Omega_k}$$
where $\Bar{w}_i^k|_{\partial\Omega_k}$ is the boundary value of $\Bar{w}_i^k$ on $\partial \Omega_k$ (constant). From (\ref{lem3.18}) and (\ref{lem3.19}) we have

$$|-\Delta \Bar{w}_i^k(\eta)|\leq\frac{C(1+|\eta|)^{-m+2\gamma+\delta}}{(1+|y_k|)^\delta}+\frac{C(1+|\eta|)^{1-m+2\gamma}}{\Lambda_k(1+|y_k|)^\delta} $$

Since we have $\Bar{w}_i^k(0)=0$, for some $i\in I$ we have

\begin{equation}\label{lem3.20}
\begin{aligned}
1&=|\Bar{w}_i^k(y_k)-\Bar{w}_i^k(0)|\\
&\leq C\int_{\Omega_k}|G(y_k,\eta)-G(0,\eta)|\left(\frac{(1+|\eta|)^{-m+2\gamma+\delta}}{(1+|y_k|)^\delta}+\frac{(1+|\eta|)^{1-m+2\gamma}}{\Lambda_k(1+|y_k|)^\delta} \right)
\end{aligned} 
\end{equation}

where the constant on the boundary is canceled out. To compute the right hand side above, we decompose the $\Omega_k$ as $\Omega_k=\Sigma_1\cup\Sigma_2\cup\Sigma_3$ where

\begin{equation*}
\begin{aligned} 
\Sigma_1&=\{\eta\in \Omega_k: |\eta|<|y/2|\}\\
\Sigma_2&=\{\eta\in \Omega_k: |y-\eta|<|y/2|\}\\
\Sigma_3&=\Omega_k\setminus(\Sigma_1\cup\Sigma_2)
\end{aligned} 
\end{equation*}
for $y\in \Omega_k$.
Using (\ref{1020e5}) we have

\begin{align*}
  \int_{\Sigma_1\cup\Sigma_2}|G(y_k,\eta)-G(0,\eta)|(1+|\eta|)^{-m+2\gamma+\delta}d\eta \\
  =O(1)(\log|y_k|)(1+|y_k|)^{(2\mu-m+\delta)_+},
\end{align*}
where
\begin{equation*}
(1+|y_k|)^{\alpha_+}=\left\{\begin{array}{ll}
  (1+|y_k|)^\alpha,&\alpha >0 \\
  \log(1+|y_k|),&\alpha =0 \\
  1,&\alpha <0 
  \end{array}\right.
\end{equation*}

$$\int_{\Sigma_3}|G(y_k,\eta)-G(0,\eta)|(1+|\eta|)^{-m+2\gamma+\delta}d\eta=O(1)(1+|y_k|)^{2\mu-m+\delta}.$$

Thus

$$\int_{\Omega_k}|G(y_k,\eta)-G(0,\eta)|\frac{(1+|\eta|)^{-m+2\gamma+\delta}}{(1+|y_k|)^\delta}d\eta=O(1)(1+|y_k|)^{2\mu-m}.$$

Similarly we can compute the other term:

$$\int_{\Omega_k}|G(y_k,\eta)-G(0,\eta)|\frac{(1+|\eta|)^{1-m+2\gamma}}{\Lambda_k(1+|y_k|)^\delta}=o(1).$$
Note that in addition to the application of estimates of $G$, we also used the assumption $\delta>1+2\mu-m$. 

\medskip

We see that the right hand side of (\ref{lem3.20}) is $o(1)$, a contradiction to the left hand side of (\ref{lem3.20}). Then Lemma \ref{three-p-1} is proved.$\Box$.

\subsection{Second order estimates}

Now we want to improve the estimates in Lemma \ref{three-p-1}. The following theorem does not distinguish $m=4\mu$ or $m<4\mu$. 

\begin{thm}
\label{2.1}
Let $\delta$ be the same as in Lemma \ref{three-p-1} there exist $C(\delta)>0$ independent of $k$ such that for $|\alpha|=0,1,$

\begin{equation}
|D^\alpha (w_i^k(y)-\Phi_i^k(y))|\leq C\bar \epsilon_k^2(1+|y|)^{2+2\mu-m-|\alpha|+\epsilon}.
\end{equation}
where
$$\Phi_i^k(y)=\bar \epsilon_k(G_{1,i}^k(r)\cos\theta+G_{2,i}^k(r)\sin\theta)$$
with
$$|G_{t,i}^k(r)|\leq Cr(1+r)^{2\mu-m+\epsilon},\quad t=1,2.$$
\end{thm}

Here we note that $\Phi$ is the projection of $w_i^k$ on $e^{i\theta}$. 

\noindent{\bf Proof:}

Here $\Phi^k(y)$ denote the projection of $v_i^k$ onto $span\{\sin\theta,\cos\theta\}$.

Taking the difference between $v^k$ and $U^k$ we have

$$
  \Delta w_i^k+\sum_{j=1}^n a_{ij}|y|^{2\gamma}\mathfrak{h}_j^k(\bar \epsilon_k y)e^{U_j^k+w_j^k}-\sum_{j=1}^n a_{ij}|y|^{2\gamma}e^{U_j^k}=0,
$$

which is
$$
  \Delta w_i^k+\sum_{j=1}^n a_{ij}|y|^{2\gamma}e^{U_j^k}(\mathfrak{h}_j^{k}(\bar \epsilon_k y)e^{w_j^k}-1)=0,
$$

We further write the equation for $w_i^k$ as

\begin{equation}
  \Delta w_i^k+\sum_{j=1}^n a_{ij}|y|^{2\gamma}e^{U_j^k}w_j^k=E_i^k,
\end{equation}

where 
\begin{align}\label{err-iw}
  E_i^k&=-\sum_{j=1}^n a_{ij}|y|^{2\gamma}e^{U_j^k}(\mathfrak{h}_j^{k}(\bar \epsilon_ky)e^{w_j^k}-1-w_j^k)\\
&=-\sum_{j=1}^n a_{ij}|y|^{2\gamma}e^{U_j^k}((\mathfrak{h}_j^{k}(\epsilon_k y)-1)+(\mathfrak{h}_j^{k}(\epsilon_k y)-1)w_j^k+O(w_j^k)^2). \nonumber\\
&=-\sum_{j=1}^n a_{ij}|y|^{2\gamma}e^{U_j^k}((\mathfrak{h}_j^{k}(\epsilon_k y)-1)+O(\varepsilon_k^2(1+|y|)^{1+\delta-m+2\gamma}) \nonumber
\end{align}
where Lemma \ref{three-p-1} is used to evaluate the last two terms.

{\bf Step one:}
We first estimate the radial part of $w_i^k$.
Let $g^{k,0}=(g_1^{k,0},\cdots,g_n^{k,0})$ be the radial part of $w_i^k$:

$$g_i^{k,0}=\frac{1}{2\pi}\int_0^{2\pi} w_i^k(r\cos\theta,r\sin\theta)d\theta$$

$g_i^k$ satisfies

\begin{equation}
L_i g^{k,0}=-\frac{\bar \varepsilon_k^2}{4}\sum_{j=1}^n a_{ij}|y|^{2\gamma}\Delta \mathfrak{h}_j^{k}(0)r^2e^{U_j^k}+O(\bar \varepsilon_k^2(1+|y|)^{1+\delta-m+2\gamma+\epsilon}).
\end{equation}

where 

$$L_i g^{k,0}=\frac{d^2}{dr^2}g_i^{k,0}+\frac{1}{r}\frac{d}{dr}g_i^{k,0}+\sum_{j=1}^n a_{ij}|y|^{2\gamma}e^{U_j^k}g_j^{k,0}$$

We claim that 
\begin{equation}\label{30}
|g_i^{k,0}|\leq C\bar \varepsilon_k^2(1+r)^{2+2\mu-m+\epsilon}, \quad 0<r<\bar \varepsilon_k^{-1}
\end{equation}
holds for some C independent of $k$. To prove (\ref{30}), we first observe that 
$$
|L_ig^{k,0}|\leq C\bar \varepsilon_k^2(1+r)^{2\mu-m}
$$
We shall contruct $f^k=(f_1^k,...,f_n^k)$ to ``replace" $L_ig^{k,0}$:
Let $f^k=(f_1^k,\cdots,f_n^k)$ be the solution of 
\begin{equation}\label{31}
\left\{\begin{array}{ll}
  \frac{d^2}{dr^2}f_i^k+\frac{1}{r}\frac{d}{dr}f_i^k=L_ig^{k,0}, \\
  \\
  f_i^k(0)=\frac{d}{dr}f_i^k(0)=0. \quad i\in I.
  \end{array}\right.
\end{equation}

The elementary estimates lead to this estimate of $f_i^k$:

\begin{equation}
|f_i^k(r)|\leq C\bar \varepsilon_k^2(1+r)^{2+2\mu-m+\epsilon}
\end{equation}

Let 
\begin{equation}\label{g1k}
g^{k,1}=g^{k,0}-f^k
\end{equation}
Clearly we have
\begin{equation}\label{33}
\left\{\begin{array}{ll}
  L_i g^{k,1}=F_i^k \\
  g^{k,1}_i(0)=\frac{d}{dr}g^{k,1}_i(0)=0. \quad i\in I.
  \end{array}\right.
\end{equation}
where 
$$F_i^k:=-\sum_{j=1}^n a_{ij}|y|^{2\gamma}e^{U_j^k}f_j^{k}=O(\bar \varepsilon_k^2)(1+r)^{4\mu-2m+\epsilon}.$$
From here we can see the purpose of $f^k$: If we can prove 
$$g_i^{k,1}(r)=O(\bar \epsilon_k^2)(1+r)^{2+2\mu-m+\epsilon}$$
the same estimate holds for $g_i^{k,0}$ because $f^k$ satisfies the same estimate. The advantage of this replacement is that now the error of $L_ig^{k,1}$ is smaller.

If $m-2\mu>\frac{1}{2}$, we have $2+4\mu-2m<1$, this is the main requirement for using Lemma \ref{three-p-1}. Employing the argument of Lemma \ref{three-p-1} we obtain

\begin{equation}\label{34}
|g^{k,1}_i(r)|\leq C\bar \varepsilon_k^2(1+r)^{2+4\mu-2m}=O(\bar \epsilon_k^2)(1+r)^{2+2\mu-m+\epsilon}. 
\end{equation}
In this case it is easy to see that we have obtained the desired estimate for $g^{k,0}$.

If $m-2\mu\le \frac{1}{2}$, we apply the same ideas by adding more correction functions to $g^{k,1}$. 
Let $\bar f^k=(\bar f_1^k,\cdots,\bar f_n^k)$ be the solution of 

$$
\left\{\begin{array}{ll}
  \frac{d^2}{dr^2}\bar f_i^k+\frac{1}{r}\frac{d}{dr}\bar f_i^k=L_ig^{k,1}, \\
  \\
  \bar f_i^k(0)=\frac{d}{dr}\bar f_i^k(0)=0. \quad i\in I.
  \end{array}\right.
$$
Then $|\bar f_i^k(r)|\le C\bar \epsilon_k^2(1+r)^{2+4\mu-2m+\epsilon}$. Let 
$$g^{k,2}_i(r)=g^{k,1}_i(r)-\bar f_i^k, \quad i\in I, $$
Then 
$$L_ig^{k,2}=O(\bar \epsilon_k^2)(1+r)^{6\mu-3m+\epsilon}. $$
If $m-2\mu>\frac 13$, we employ the method of Lemma \ref{three-p-1} to obtain 
$$|g^{k,2}_i(r)|\le C\bar \epsilon_k^2(1+r)^{2+6\mu-3m}=O(\bar \epsilon_k^2)(1+r)^{2+2\mu-m+\epsilon}, $$
and the proof is complete if $m-2\mu>\frac 13$. Obviously if $m-2\mu\le \frac{1}{3}$, such a correction can be done finite times until (\ref{30}) is eventually established.

\medskip

\noindent{\bf Step 2: Projection on $\sin\theta$ and $\cos\theta$:}

\medskip

In this step we consider the projection of $w_i^k$ over $\cos\theta$ and $\sin\theta$ respectively:

$$\bar \varepsilon_k G_{1,i}^k(r)=\frac{1}{2\pi}\int_0^{2\pi} w_i^k(r,\theta)\cos\theta d\theta,\quad \bar \varepsilon_k G_{2,i}^k(r)=\frac{1}{2\pi}\int_0^{2\pi} w_i^k(r,\theta)\sin\theta d\theta$$

Clearly, $G_{1,i}^k(r)$ and $G_{t,i}^k(r)$ solve the following linear systems for $t=1,2$ and $r\in (0,\bar \varepsilon_k^{-1})$:

\begin{equation}\label{35}
\begin{aligned}
  &(\frac{d^2}{dr^2}+\frac{1}{r}\frac{d}{dr}-\frac{1}{r^2})G_{t,i}^k+\sum_ja_{ij}|y|^{2\gamma}e^{U_j^k}G_{t,j}^k \\
  &=-\sum_j a_{ij}|y|^{2\gamma} \partial_t \mathfrak{h}_j^{k}(0)re^{U_j^k}+O(\varepsilon_k)(1+r)^{1+\delta+2\gamma-m+\epsilon}
  \end{aligned}
\end{equation}
where $\delta=1-\epsilon$ if $m-2\mu\ge 1$ and $\delta=1+2\mu-m+\epsilon$ if $m-2\mu\le 1$.  Thus 
$$(\frac{d^2}{dr^2}+\frac{1}{r}\frac{d}{dr}-\frac{1}{r^2})G_{t,i}^k+\sum_ja_{ij}|y|^{2\gamma}e^{U_j^k}G_{t,j}^k=O(r)(1+r)^{2\gamma-m+\epsilon}. $$
Let 

\begin{equation}
\Phi_i^k=\bar \varepsilon_k G_{1,i}^k(r)\cos\theta+\bar \varepsilon_k G_{2,i}^k(r)\sin\theta
\end{equation}

Then $\Phi^k$ solves

\begin{equation}
\begin{aligned}
  &\Delta \Phi_i^k+\sum_ja_{ij}|y|^{2\gamma}e^{U_j^k}\Phi_i^k \\
  &=-\bar \varepsilon_k\sum_j a_{ij}|y|^{2\gamma} (\partial_1 \mathfrak{h}_j^{k}(0)y_1+\partial_2 \mathfrak{h}_j^{k}(0)y_2)e^{U_j^k}+O(\bar \varepsilon_k^2)(1+r)^{1+\delta+2\gamma-m+\epsilon}
  \end{aligned}
\end{equation}

By Lemma 2.2 we have,
$$|G_{1,i}^k(r)|+|G_{2,i}^k(r)|\leq C(1+r)^\delta$$

Then we can rewrite (\ref{35}) as

\begin{equation}
(\frac{d^2}{dr^2}+\frac{1}{r}\frac{d}{dr}-\frac{1}{r^2})G_{t,i}^k=h(r),\quad r\in (0,\bar \varepsilon_k^{-1})
\end{equation}

where 
$$|h(r)|\leq C(1+r)^{1-m+2\gamma+\epsilon}.$$

By standard ODE theory
\begin{equation}\label{39}
G_{1,i}^k=c_{1k}r+\frac{c_{2k}}{r}-\frac{r}{2}\int_r^\infty h(s)ds-\frac{r^{-1}}{2}\int_0^r s^2h(s)ds
\end{equation}

Since $G_{1,i}$ is bounded near $0$, $c_{2k}=0$.
Using $G_{1,i}^k(\bar \varepsilon_k^{-1})=O(\bar \varepsilon_k^{-\delta})$, we have

$$|c_{1k}|\leq C\bar \varepsilon_k^{1-\delta}$$

Then from $(\ref{39})$

\begin{equation}\label{40}
|G_{t,i}^k(r)|\leq Cr(1+r)^{2\mu-m+\epsilon}+c\bar \epsilon_k^{1-\delta}r,\quad t=1,2.
\end{equation}
If $m-2\mu\le 1$, the last term can be ignored since it is not greater than the first term on the right. When $m-2\mu>1$, $\delta=1-\epsilon$. In this situation we need to keep this term at this moment. 

\medskip

\noindent{\bf Step three: Projection onto higher notes.}

\medskip

Let $g^{k,l}=(g_1^{k,l},\cdots,g_n^{k,l})$ be the projection of $w^{1,k}$ on $\sin l\theta$. 

First we prove a uniqueness lemma for global solutions: 
\begin{lem}\label{uniq-high-nodes}
Let $g=(g_1,...,g_n)$ satisfy
$$g_i''(r)+\frac 1r g_i'(r)-\frac{l^2}{r^2}g_i+\sum_j a_{ij} r^{2\gamma} e^{U_j}g_j=0,\quad 0<r<\infty, 
\quad i=1,..,n$$
and $g_i(r)=o(r)$ for $r$ small, $g_i(r)=O(r^{2+2\mu-m})$ for $r$ large, then $g_i\equiv 0$ if $l\ge 2$. 
\end{lem}

\noindent{\bf Proof of Lemma \ref{uniq-high-nodes}:} Treating $\sum_j a_{ij}r^{2\gamma}e^{U_j}g_j$ as an error term, standard ode theory gives
\begin{align*}
g_i(r)&=c_1r^l+c_2r^{-l}\\
&+r^l\int_{\infty}^r \frac{\sum_j a_{ij}s^{2\gamma-l}e^{U_j(s)}g_j(s)}{(-2l)/s}ds
-r^{-l}\int_0^r\frac{s^l\sum_j a_{ij}s^{2\gamma}e^{U_j}g_j(s)}{(-2l)/s}ds. 
\end{align*}
From the bound near $0$, $c_1=0$, from the bound at infinity $c_2=0$. Thus standard evaluation gives 
$$|g_i(r)|\le C(1+r)^{2+4\mu-2m}$$ for some $C>0$ and all $r$. Thus if $m-2\mu\ge \frac 12$, $g_i(r)\cos(l\theta)$
would be a bounded solution of 
$$\Delta \phi_i+\sum_j a_{ij}|x|^{2\gamma}e^{U_j}\phi_j=0,\quad \mbox{in}\quad \mathbb R^2, \quad i=1,..,n. $$
Since we have also $\phi_i(0)=0$ for all $i$,
 Proposition \ref{uniqlin} gives $\phi_i\equiv 0$, which is $g\equiv 0$. If $m-2\mu<\frac 12$, we use the new bound of $g_i$ to improve it to 
 $$|g_i(r)|\le C(1+r)^{2+6\mu-3m}.$$
 Obviously it takes finite steps to prove that $g_i$ is bounded. Lemma \ref{uniq-high-nodes} is established. $\Box$

\medskip

Next for each fixed $l\ge 2$ we prove an estimate for $g^{k,l}$:

\begin{lem}\label{est-high} For each $l\ge 2$ fixed, there exists $C(l)>0$ such that 
$$|g^{k,l}_i(r)|\le C(l)\bar \epsilon_k^2(1+r)^{2+2\mu-m+\epsilon}. $$
\end{lem}

\noindent{\bf Proof of Lemma \ref{est-high}}: 

If the estimate is false, we would have 
$$\Lambda_k=\max_i\max_r\frac{|g_i^{k,l}(r)|}{\bar \epsilon_k^2(1+r)^{2+2\mu-m+\epsilon}}\to \infty. $$
Suppose $\Lambda_k$ is attained at $r_k$. Then we set $w_i^k$ as
$$w_i^k(r)=\frac{g_i^{k,l}(r)}{\Lambda_k\bar \epsilon_k^2(1+r_k)^{2+2\mu-m+\epsilon}}. $$
From the definition of $\Lambda_k$ we see that
\begin{equation}\label{w-ik-bound}
|w_i^k(r)|\le \frac{(1+r)^{2+2\mu-m+\epsilon}}{(1+r_k)^{2+2\mu-m+\epsilon}}.
\end{equation}
The equation of $w_i^k$ is
$$\frac{d^2}{dr^2}w_i^k+\frac 1r \frac{d}{dr}w_i^k-\frac{l^2}{r^2}w_i^k+\sum_j a_{ij}r^{2\gamma}e^{U_j}w_j^k=
\frac{o(1)r^{2\mu}(1+r)^{-m}}{(1+r_k)^{2+2\mu-m}}. $$
First we rule out the case that $r_k$ converges somewhere, because this case would lead to a solution $w=(w_1,...,w_n)$ of 
$$w_i''(r)+\frac 1rw_i'(r)+\sum_j a_{ij} r^{2\gamma}e^{U_j}w_j=0,\quad 0<r<\infty, \quad i=1,..,n $$
with $w_i$ bounded. But according to Proposition \ref{uniqlin} $w_i\equiv 0$, impossible to have $\Lambda_k$ attained at $r_k$. 

Next we rule out the case $r_k\to \infty$. On one hand we have $\max_i|w_i^k(r_k)|=1$. On the other hand, the evaluation of $w_i^k(r_k)$ gives
$$w_i(r)=c_1r^l+c_2r^{-l}+r^l\int_{\infty}^r\frac{f(s)s^{1-l}}{2l}ds-r^{-l}\int_0^r\frac{s^{l+1}f(s)}{2l}ds $$
where 
$$f(r)=-\sum_j a_{ij}r^{2\gamma}e^{U_j(r)}w_j^k(r)+o(1)\frac{r^{2\mu}(1+r)^{-m+\epsilon}}{(1+r_k)^{2+2\mu-m+\epsilon}}.$$
First we observe that $c_2=0$, $c_1=O(\bar \epsilon_k^{l+1-\delta})$ for 
$\delta=1-\epsilon$ if $m-2\mu>1$ and $\delta=1+2\mu-m+\epsilon$ if $m-2\mu\le 1$. Using the bound in (\ref{w-ik-bound})
we see that $w_i^k(r_k)=o(1)$, which is a contradiction to $\max_i|w_i^k(r_k)|=1$. Lemma \ref{est-high} is established. 

\medskip

Next we obtain a uniform estimate for all projections in high nodes. We claim that there exists $C>0$ independent of $k$ and $l\ge 2$ such that
\begin{equation}\label{unif-l}
|g_i^{k,l}(r)|\le C\bar \epsilon_k^2(1+r)^{2+m-2\mu+\epsilon}+C\bar \epsilon_k^{l+1-\delta}r^l,\quad 0<r<\tau \ \bar \epsilon_k^{-1}. 
\end{equation}
Obviously we only need to consider large $l$. 
Let 
$$F= (1+r)^{2+m-2\mu+\epsilon}+\bar \epsilon_k^{l-1-\delta}r^l$$
and $g$ be defined as
$$g_i(r)=\frac{r^2}4\int_r^{\infty}Q\bar \epsilon_k^2s^{-1}F(s)ds+
\frac 14r^{-2}\int_0^rQ\bar \epsilon_k^2F(s)s^3ds, $$
for all $i=1,..,n$, then all $g_i$ satisfies
$$g_i''(r)+\frac 1rg_i'(r)-\frac{4}{r^2}g_i(r)=-Q\bar \epsilon_k^2F. $$
By choosing $Q$ large we have 
$L_ig_i\le L_ig^{k,l}$ and we can write the equation of $g_i$ as
$$g_i''+\frac 1rg_i'(r)-\frac{l^2}{r^2}g_i=-Q\bar \epsilon_k^2F(r)
+\frac{4-l^2}{r^2}g_i<-Q\bar \epsilon_k^2F(r).$$
We also observe that $g_i(\bar\epsilon_k^{-1})>g^{k,l}_i(\bar\epsilon_k^{-1})$. Thus by setting $h_i=g_i-g^{k,l}_i$, we have 
$$h_i''(r)+\frac 1rh_i'(r)-\frac{l^2}{r^2}h_i+\sum_j a_{ij}r^{2\gamma}e^{U_j}h_j\le 0
\quad 0<r<\bar\epsilon_k^{-1}. $$
Here we claim that $h_i\ge 0$. If not, without loss of generality $\min_i\min_r h_i(r)$ is attained at $h_1(r_k)<0$.  Looking at the equation for $h_1^k$:
$$h_1''(r_k)+\frac{1}{r}h_1'(r_k)+\sum_{j=1}^na_{ij}r^{2\gamma}e^{U_j}h_j(r_k)-\frac{l^2}{r^2}h_1(r_k)\le 0. $$
Since $h_i(r_k)\ge h_1^k(r_k)$. For $l_0$ large and $l\ge l_0$, the left hand side is positive. A contradiction. Thus (\ref{unif-l}) is established.

\medskip

Next we obtain for $l\ge 2$ that 
\begin{equation}\label{depend-l}
|g_i^{k,l}(r)|\le \frac{C}{l^2}\bar \epsilon_k^2(1+r)^{2+2\mu-m+\epsilon}
+c\bar \epsilon_k^{l-\delta}r^l, \quad 0<r<\bar\epsilon_k^{-1}. 
\end{equation}
From the expression of $w_i^k$ we have
$$\frac{d^2}{dr^2}g_i^{k,l}+\frac 1r\frac{d}{dr}g_i^{k,l}+\sum_j a_{ij}r^{2\gamma}e^{U_j}g_j^{k,l}-\frac{l^2}{r^2}g_i^{k,l}=f_1 $$
where
$$|f_1|\le C\bar \epsilon_k^2(1+r)^{2\mu-m+\epsilon}. $$
Setting 
$$f_2=f_1-\sum_j a_{ij}r^{2\gamma}e^{U_j}g_j^{k,l}, $$
then the estimate of $g_i^{k,l}$ we have 
$$f_2=O(\bar \epsilon_k^2)(1+r)^{2\mu-m+\epsilon}. $$
we write the equation for $g_i^{k,l}$ as 
$$\frac{d^2}{dr^2}g_i^{k,l}+\frac 1r\frac{d}{dr}g_i^{k,l}-\frac{l^2}{r^2}g_i^{k,l}=f_2 $$
Thus we have
$$g_i^{k,l}=c_1r^l+r^l\int_{\infty}^r\frac{(-f_2)s^{-l}}{(-2l)/s}ds
+r^{-l}\int_0^r\frac{s^lf_2(s)}{(-2l)/s}ds
$$
Then it is easy to see that the last two terms are of the order $O(\bar \epsilon_k^2)(1+r)^{2+2\mu-m+\epsilon}/l^2$. To obtain the order of $c_1$, we use the boundary value to finish the proof of (\ref{depend-l}). 
On $r=\tau \bar \epsilon_k^{-1}$ the crude estimate in (\ref{11mar9e2}) gives $c_1=O(\bar \epsilon_k^{l+1-\delta})$. After using (\ref{unif-l}) in the evaluation above, (\ref{depend-l}) can be obtained by elementary estimates.

Taking the sum of $g_i^{k,l}$ we have obtained the following estimate of $w^k$:
\begin{equation}\label{tem-est}
|(w_i^k-\Phi_i^k)(x)|\le C\bar \epsilon_k^2(1+|y|)^{2+2\mu-m+\epsilon}+C\bar \epsilon_k^{2-\delta}|y|^{2}. 
\end{equation}
where $\delta=1-\epsilon$ if $m-2\mu>1$ and $\delta=1+2\mu-m+\epsilon$ if $m-2\mu\le 1$. Now we improve the estimate of $w_i^k$ by considering its spherical average: 
$$\bar w_i^k(r)=\frac 1{2\pi r}\int_{B_r}w_i^k. $$
$$\frac{d}{dr}\bar w_i^k(r)=\frac{1}{2\pi r}\int_{B_r}\Delta w_i^k. $$
Using (\ref{tem-est}) in the evaluation of the right hand side we have
$$\frac{d}{dr}\bar w_i^k(r)=O(\bar \epsilon_k^2)/r. $$
Thus the value of $w_i^k$ on $\partial B(0,\tau)$ is $O(\bar \epsilon_k^2)\log \frac{1}{\bar \epsilon_k}$. 
Here we recall that $w_i^k$ on the outside boundary is a constant. Using this new information in the proof we see that the $O(\epsilon_k^{2-\delta}r^2)$ can be removed. 
The statement of Theorem \ref{2.1} is established for $\alpha=0$.

Once the statement for $|\alpha |=0$ is obtained, the statement for $|\alpha |=1$ follows from standard bootstrap argument. Obviously we only need to discuss the regularity around the origin. First $\gamma>-1$ implies that $w_k\in W^{2,1+\epsilon}$ for some $\epsilon>0$, thus $w_k\in C^{\tau}$ for some $\tau>0$. Using this improve regularity and $w_k(0)=0$ we observe that $w_k\in W^{2,1+\epsilon_1}$ for a larger $\epsilon_1>0$, which results in more smoothness of $w_k$. After finite steps we have $w_k\in W^{2,2+\epsilon_2}$ for some $\epsilon_2>0$, which leads to the desired first order estimates. 
Theorem \ref{2.1} is established. $\Box$

\end{document}